\input amstex
\input epsf
\documentstyle{amsppt}
\NoRunningHeads
\NoBlackBoxes
\document

\define\ci{\circledcirc}
\def\O{\Cal O}

\def\Ua{U_q(\tilde\g)}
\def\U2{{\Ua}_2}
\def\g{\frak g}

\def\Z{\Bbb Z}
\def\C{\Bbb C}

\def\l{\lambda}

\def\<{\langle}
\def\>{\rangle}
\def\o{\otimes}

\def\End{\text{End}}

\topmatter
\title  Set-theoretical solutions to the quantum Yang-Baxter equation
\endtitle

\author Pavel Etingof, Travis Schedler,  Alexandre Soloviev
\endauthor

\address
\newline
P.~E. : Department of Mathematics, Harvard University, Cambridge, MA
02138
\newline
T.~S. : 2312 Harvard Yard Mail Center, Cambridge, MA 02138
\newline
A.~S. : Department of Mathematics, MIT, Cambridge, MA 02139
\newline
\newline
Email addresses:
\newline
P.~E. : $etingof\@math.harvard.edu$
\newline
T.~S. : $schedler\@fas.harvard.edu$
\newline
A.~S. : $sashas\@math.mit.edu$
\endaddress

\abstract
In the paper \cite{Dr}, V.Drinfeld formulated a number of problems in
quantum group theory.  In particular, he suggested to consider
``set-theoretical'' solutions of the quantum Yang-Baxter equation,
i.e. solutions given by a permutation $R$ of the set $X\times X$,
where $X$ is a fixed set. In this paper we study such solutions, which
in addition satisfy the unitarity and nondegeneracy conditions. We
discuss the geometric and algebraic interpretations of such solutions,
introduce several constructions of them, and give their classification
in group-theoretic terms.
\endabstract
\endtopmatter
\head 0. Introduction\endhead

The quantum Yang-Baxter equation is one of the basic equations in
mathematical physics, which lies in the foundation of the theory of
quantum groups. This equation involves a linear operator $R: V\o V\to
V\o V$, where $V$ is a vector space, and has the form
$$
R^{12}R^{13}R^{23}=R^{23}R^{13}R^{12}\text{ in }\text{End}(V\o V\o V),
$$
where $R^{ij}$ means $R$ acting in the i-th and j-th components. 

In the last 15 years, many solutions of this equation were found and
the related algebraic structures (Hopf algebras) have been intensively
studied. However, these solutions were usually ``deformations'' of the
identity solution. On the other hand, it is interesting to study
solutions which are not obtained in this way. In \cite{Dr}, Drinfeld
suggested to study the simplest class of such solutions -- the so
called set-theoretical solutions.  By definition, a set-theoretical
solution is a solution for which $V$ is a vector space spanned by a
set $X$, and $R$ is the linear operator induced by a mapping $X\times
X\to X\times X$.
  
In this paper we study set-theoretical solutions of the quantum
Yang-Baxter equation, satisfying additional conditions: invertibility,
unitarity, and nondegeneracy.  They turn out to have many beautiful
properties.  We discuss the geometric and algebraic interpretations of
such solutions, introduce several constructions of them, and give
their classification in terms of group theory.

The brief content of the paper is as follows. 

Chapter 1 contains the background material.  In Section 1.1 we give
the main definitions and the simplest examples. We introduce the
notion of a nondegenerate symmetric set, which is a set $X$ with an
invertible mapping $R:X^2\to X^2$ satisfying the quantum Yang-Baxter
equation and the nondegeneracy and unitarity conditions.  We explain
that if $X$ is a nondegenerate symmetric set then the set $X^n$ has a
natural action of the symmetric group $S_n$, called the twisted
action, which is, in general, different from the usual action by
permutations.  In Section 1.2 we show that any nondegenerate symmetric
set defines a coloring rule for collections of closed smooth curves in
the plane, under which the number of colorings depends only on the
number of curves involved, and not at the pattern of their
intersections.  This gives a topological interpretation of the notion
of a nondegenerate symmetric set. In Section 1.3 we show that the
twisted action of $S_n$ on $X^n$ for a nondegenerate symmetric set is
conjugate to the action by permutations.

Chapter 2 introduces and studies the main algebraic structure
associated to a nondegenerate symmetric set $X$ -- its structure group
$G_X$. In Section 2.1 we show that $G_X$ has two natural actions on
$X$, which are conjugate to each other. In Sections 2.2,2.3 we show
that the group $G_X$ is naturally a subgroup of $Aut(X)\ltimes \Z^X$,
such that the 1-cocycle defined by the projection $G_X\to \Z^X$ is
bijective. Using this result, in Section 2.4 we show that
nondegenerate symmetric sets, up to isomorphism, are in 1-1
correspondence with quadruples $(G,X,\rho,\pi)$, where $G$ is a group,
$X$ is a set, $\rho$ a left action of $G$ on $X$, and $\pi$ a
bijective 1-cocycle of $G$ with coefficients in $\Z^X$. In Sections
2.5-2.6 we show that there exists a unique, up to isomorphism,
indecomposable nondegenerate symmetric set of order $p$, where $p$ is
a prime -- $X=\Z/p\Z$, $R(x,y)=(x+1,y-1)$. In Section 2.7 we prove
solvability of the structure group.  In Section 2.8 we apply the
notion of the structure group of a nondegenerate symmetric set.  to
the study of decomposable nondegenerate symmetric sets.  Finally, in
Sections 2.9-2.10 we study the quantum algebras associated to a
nondegenerate symmetric set by the Faddeev-Reshetikhin-
Takhtajan-Sklyanin construction.

Chapter 3 introduces the main constructions of nondegenerate symmetric sets --
linear, affine, multipermutation solutions, twisted unions, 
generalized twisted unions. In this chapter we classify 
such solutions, and study their properties. At the end we give 
the results of a computer calculation, which found all 
nondegenerate symmetric sets $X$ with $|X|\le 8$. 
%It turned out that all 
%solutions found belong to one of the types defined in Chapter 3. 
%We don't expect, however, that this will necessarily be the case for larger 
%$|X|$ or infinite X.  

In Chapter 4 we consider power series solutions 
of the Yang-Baxter equation, which are a generalization 
of linear solutions. We show that a power series solution 
with a generic linear part is equivalent to a linear solution.

In the appendix we introduce the notion of a $T$-structure
on an abelian group $A$, which is motivated by the definition of the map $T$
in Proposition 2.2. We discuss the connection of $T$-structures with bijective 
1-cocycles, in particular in the case of a cyclic group $A$. 

\head Acknowledgments \endhead

We are grateful to Noam Elkies for useful discussions.
The work of Pavel Etingof was supported by an NSF grant.
Travis Schedler thanks MIT for hospitality. 
  
\head 1. Braided and symmetric sets\endhead

\subhead 1.1. Definitions\endsubhead

Let $X$ be a nonempty set, and $S:X\times X\to X\times X$ be a bijection.
We will denote the components of $S$ by $S_1$ and $S_2$
(i.e. $S(x_1,x_2)=(S_1(x_1,x_2),S_2(x_1,x_2))$); they are binary 
operations on $X$. For positive integers $i<n$ 
let the map $S^{ii+1}:X^n\to X^n$ be defined by 
$S^{ii+1}=id_{X^{i-1}}\times S\times id_{X^{n-i-1}}$. 

\proclaim{Definition 1.1} 

(i) A pair $(X,S)$ is called nondegenerate if 
the maps $X\to X$ defined by 
$x\to S_2(x,y)$ and $x\to S_1(z,x)$ are bijections for any fixed $y,z\in X$. 
 
(ii) A pair $(X,S)$ is said to be a braided set if $S$ satisfies
the braid relation 
$$
S^{12}S^{23}S^{12}=S^{23}S^{12}S^{23}.\tag 1.1
$$

(iii) A pair $(X,S)$ is called involutive if 
$$
S^2=id_X\tag 1.2
$$
A braided set $(X,S)$ which is involutive is called a symmetric set. 

(iv) Pairs $(X,S)$ and $(X',S')$ are said to be isomorphic 
if there exists a bijection $\phi:X\to X'$ which maps $S$ to $S'$.  
\endproclaim

The main objects of study in this paper are nondegenerate symmetric sets. 
Our main goal is to learn to construct them and to understand 
their properties. For brevity, nondegenerate symmetric sets will often be 
called ``solutions'' (meaning nondegenerate solutions of equations
(1.1),(1.2)). 

{\bf Examples.} 1. Let $X$ be any set, and $S(x,y)=(y,x)$. 
Then $(X,S)$ is a nondegenerate symmetric set. 
It is called ``the trivial solution''. 

2. (Lyubashenko, see \cite{Dr}) 
Let $X$ be any set, and $S(x,y)=(f(y),g(x))$, where $f,g:X\to X$. 
Then: $(X,S)$ is nondegenerate iff $f,g$ are bijective; 
$(X,S)$ is braided iff $fg=gf$; $(X,S)$ is involutive iff $g=f^{-1}$
(in this case it is also braided, i.e. symmetric).
In the last case $(X,S)$ is called ``a permutation solution''. 
If $f$ is a cyclic permutation, we will say that $(X,S)$ is a cyclic 
permutation solution. It is clear that two permutation solutions 
are isomorphic if and only if the corresponding permutations are conjugate. 

3. Let $(X,S_X),(Y,S_Y)$ be two solutions. 
Then $(X\times Y,S_X\times S_Y)$ is a solution, which is called the Cartesian
product of $X$ and $Y$.
 
Recall that the braid group $B_n$ is generated by elements 
$b_i, 1\le i\le n-1$, with defining relations 
$$
b_ib_j=b_jb_i,\ |i-j|>1;\ b_ib_{i+1}b_i=b_{i+1}b_ib_{i+1}, \tag 1.3
$$
and that the symmetric group $S_n$ is the quotient of $B_n$ by the relations 
$b_i^2=1$. Therefore, we have the following obvious proposition. 

\proclaim{Proposition 1.1} (i) The assignment 
$b_i\to S^{ii+1}$ extends to an action of $B_n$ on $X^n$ if and only if 
$(X,S)$ is a braided set. 

(ii) The assignment 
$b_i\to S^{ii+1}$ extends to an action of $S_n$ on $X^n$ if and only if 
$(X,S)$ is a symmetric set. 
\endproclaim

This proposition explains our terminology. 

\proclaim{Definition 1.2} The action of $B_n$ (or $S_n$) on $X_n$ defined 
by Proposition 1.1 will be called the twisted action. 
\endproclaim

Let $\sigma:X\times X\to X\times X$ be the permutation map, defined by
$\sigma(x,y)=(y,x)$. Let $R=\sigma \circ S$. The map $R$ is called 
the R-matrix corresponding to $S$. We have the following obvious proposition: 

\proclaim{Proposition 1.2} (i) $(X,S)$ is a braided set if and only if
$R$ satisfies the quantum Yang-Baxter equation 
$$
R^{12}R^{13}R^{23}=R^{23}R^{13}R^{12},\tag 1.4
$$
and is a symmetric set if and only if in addition to (1.4) $R$ satisfies 
the unitarity condition
$$
R^{21}R=1. \tag 1.5
$$
\endproclaim

A nice corollary of the properties of nondegenerate involutive pairs
$(X,S)$ is the following crossing symmetry property. 

Let $X_1,...,X_n$ be sets, $Y=X_1\times ...\times X_n$, and 
let $Q:Y\to Y$ 
be a map, $Q=(Q_1,...,Q_n)$. Suppose that 
the function $Q_i(x_1,...,x_n)$, regarded as a function of $x_i$ 
when other $x_j$ are fixed, is a bijection $X_i\to X_i$. In this case, define 
$Q^{t_i}:Y\to Y$ (the transposition of $Q$ in the i-th component)
by the condition:
if $Q(x_1,...,x_n)=(y_1,...,y_n)$ then 
$Q^{t_i}(x_1,...x_{i-1},y_i,x_{i+1},...,x_n)=
(y_1,...,y_{i-1},x_i,y_{i+1},...,y_n)$. 

\proclaim{Proposition 1.3} If $(X,S)$ is nondegenerate and 
involutive, and $R=\sigma S$, then $R^{t_i}$ is defined for $i=1,2$, and
$R$ has the crossing symmetry property
$$
R^{t_1}(R^{21})^{t_1}=R^{t_2}(R^{21})^{t_2}=id_{X^2}.\tag 1.6 
$$
In particular, $R^{t_i}$ are bijections. 
\endproclaim

\demo{Proof}
The statement of the proposition is equivalent to the equality
$$
|\{(k,l)\in X^2:S(l,j)=(k,i),S(k,i')=(l,j')\}|=\delta_{ii'}\delta_{jj'}.\tag 1.7
$$
Let us check this equality.   
if $S(l,j)=(k,i),S(k,i')=(l,j')$, then by $S^2=1$ we have
$S(k,i)=(l,j)$, which by nondegeneracy implies $i=i',j=j'$. 
Conversely, if $i=i',j=j'$, then there exist unique 
$k,l$ such that conditions (1.7) are satisfied. 
$\square$\enddemo  

{\bf Remark.} The operation $t_i$ can be defined for all mappings 
$Q$, not necessarily such that $Q_i$ is invertible as a function of $x_i$. 
To do this, we should regard $Q$ not as a map of $X_1\times...\times X_n$ to 
itself, but as a linear operator on the vector space $V_1\o...\o V_n$, 
where $V_i$ is the vector space spanned by $X_i$
(for simplicity we assume that $X_i$ are finite). In this case, 
$Q^{t_i}$ can be defined to be the endomorphism of 
$V_1\o...\o V_i^*\o...\o V_n$, obtained by dualizing the i-th component of 
$Q$. Under this definition, the crossing symmetry 
equations (1.6) make sense for any map $R$, 
and it is easy to show that if $R$ satisfies the unitarity condition, then 
crossing symmetry is equivalent to nondegeneracy. 

\subhead 1.2. Colorings of flat links\endsubhead

Nondegenerate symmetric sets turn out to have a nice geometric interpretation, 
which is given below. This interpretation is not new, but is a very simple
special case of the theory of quantum invariants of links. This theory is 
described in detail in several textbooks, e.g. \cite{Tu}.  

By a nondegenerate smooth curve in the plane we mean a parameterized 
curve $\gamma(t)=(x(t),y(t))$ such that the functions $x,y$ are smooth, and 
their derivatives are never simultaneously zero.
A nondegenerate smooth curve has a canonical orientation, defined by 
the direction of the tangent vector $\gamma'$.  

By a flat link we mean a finite 
collection of closed nondegenerate smooth curves
in the plane.

It is clear that the only singularities of a generic flat link 
are simple crossings. Thus, from combinatorial point of view, 
a generic flat link is the same thing as an oriented flat graph, 
whose vertices are all 4-valent, and have the form

\epsfbox{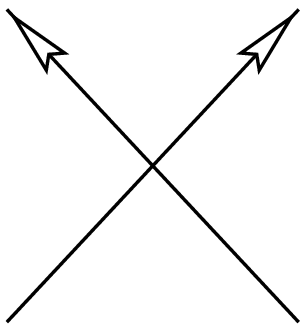}

This graph is allowed to have closed edges, without any vertices on them.

Let $X$ be a set, and $S:X^2\to X^2$ a mapping. 

\proclaim{Definition 1.3} An $X$-coloring of a generic flat link $L$ is 
an assignment to every edge of the graph $L$ of a color (an element of 
$X$), such that for any vertex of $L$ of the form 

\epsfbox{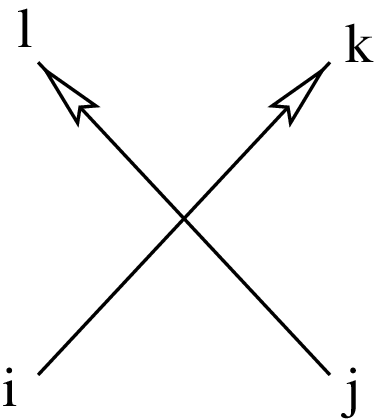}

one has $S(i,j)=(l,k)$. 
\endproclaim

For a general map $S$, it is not obvious why at least one coloring of $L$ 
exists. However, we have the following proposition. 

\proclaim{Proposition 1.4} If $(X,S)$ is a finite 
nondegenerate symmetric set, 
then the number of $X$-colorings of $L$ equals to $|X|^{n(L)}$, 
where $|X|$ is the size of $X$, and $n(L)$ the number of components in $L$. 
\endproclaim

\demo{Proof} If $L$ consists of $n(L)$ non-intersecting closed simple curves, 
then the result is clear (any component can have any color).
However, it is well known that any generic flat link can be brought 
to this form by using a sequence of the following 
Reidemeister moves:

\epsfbox{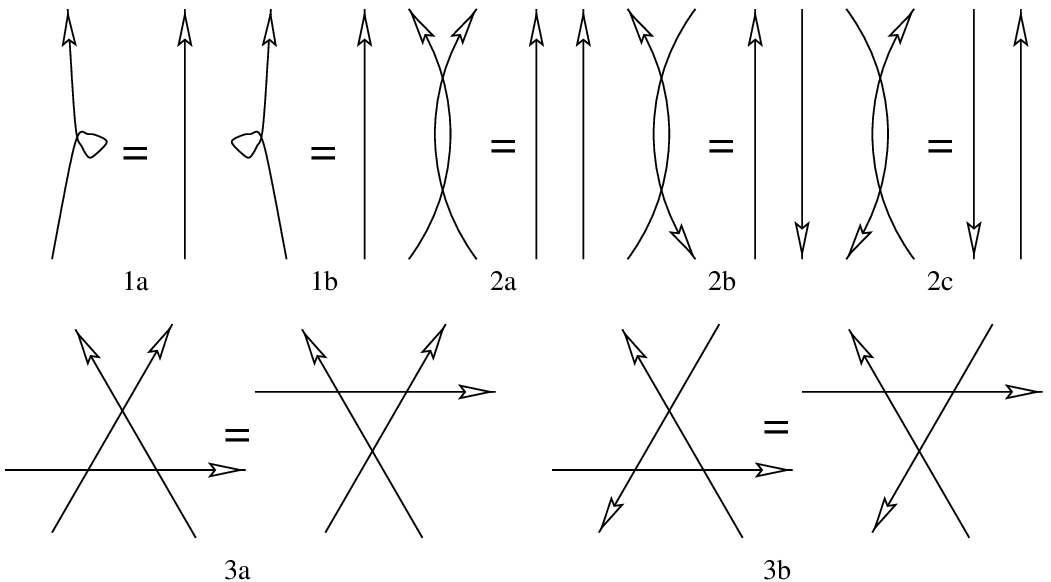}

Thus, the only thing we have to prove is that the number of colorings 
does not change when either of these moves is applied.  

The invariance of the number of colorings under move 1a
follows from the statement 
$$
|\{k\in X:S(i,k)=(j,k)\}|=\delta_{ij}.\tag 1.8
$$
This statement follows from nondegeneracy and involutivity. 
Indeed, suppose $S(i,k)=(j,k)$, then $S^2=1$ implies 
$S(j,k)=(i,k)$, which by nondegeneracy implies $i=j$. 
Thus, $S(i,k)=(i,k)$. By nondegeneracy, $k$ is unique if it exists. 
So it remains to show that for any $i$ there exists $k$ such that
$S(i,k)=(i,k)$. To do this, let $k$ be defined by $S_1(i,k)=i$. 
By nondegeneracy, such a $k$ exists. Then $S(i,k)=(i,k')$. Then 
$S(i,k')=(i,k)$ (since $S^2=1$), and thus by nondegeneracy $k=k'$. 
Move 1b is treated analogously.

The invariance under move 2a follows directly from the involutivity.
The invariance under moves 2b and 2c follows from the 
crossing symmetry. 

The invariance under move 3a is exactly the braid relation. 
The invariance under move 3b follows from crossing symmetry   
and the braid relation. Indeed, we have to check the equality
$$
(R^{21})^{t_2}R^{13}(R^{32})^{t_2}=(R^{32})^{t_2}R^{13}(R^{21})^{t_2}.\tag 1.9
$$
This can be rewritten as
$$
((R^{32})^{t_2})^{-1}(R^{21})^{t_2}R^{13}=R^{13}(R^{21})^{t_2}
((R^{32})^{t_2})^{-1}.\tag 1.10
$$
Using crossing symmetry, we have $(R^{t_i})^{-1}=(R^{21})^{t_i}$, so 
(1.10) reduces to 
$$
(R^{23})^{t_2}(R^{21})^{t_2}R^{13}=R^{13}(R^{21})^{t_2}
(R^{23})^{t_2}.\tag 1.11
$$
Transposing the second component on both sides, we get
$$
R^{21}R^{23}R^{13}=R^{13}R^{23}R^{21}, \tag 1.12
$$
which is the Yang-Baxter equation with 1 and 2 permuted. 
The proposition is proved. 
$\square$\enddemo

Proposition 1.4 has a rather trivial, but curious application. 
Suppose we have a system $L$ of closed nondegenerate smooth curves with simple 
intersections in a space of dimension $>2$, with orientations
at intersection points  (i.e. it is 
agreed which incoming edge at each vertex is the left incoming edge).
For such a system we can compute the number of colorings
as explained above. 

\proclaim{Corollary 1.5} If the number of $X$-colorings of $L$ is not equal 
to $|X|^{n(L)}$, then $L$ cannot be put on the plane 
preserving orientations at vertices, without additional self-intersections. 
\endproclaim

{\bf Example.} Consider the two-component graph $L$ of the form 

\epsfbox{labeled-si.eps}
%\vskip 1in

It is easy to see that the number of colorings of this graph 
equals to the number of fixed points of $R$ on $X^2$. So, if $R\ne id$, 
the number of colorings is less than $|X|^2$. Thus, $L$ cannot be put
on the plane without additional self-intersections (of course, this
is obvious from the picture). 

\subhead 1.3. The isomorphism of the two $S_n$-actions\endsubhead

Let $(X,S)$ be a set with a mapping.
Introduce the notation 
$$
S(x,y)=(g_x(y),f_y(x)).\tag 1.13
$$ 

\proclaim{Proposition 1.6} If $S$ is involutive then 
$$
f_{f_y(x)}g_x(y)=y.\tag 1.14
$$
If $(X,S)$ is a braided set then 
$$
f_yf_x=f_zf_t, \  
g_xg_y=g_tg_z\text{ when }S(x,y)=(t,z).\tag 1.15
$$
\endproclaim

\demo{Proof} Straightforward.
$\square$\enddemo

Recall that in the previous chapter for any symmetric set $(X,S)$ we defined 
the twisted action of $S_n$ on $X^n$. We will need the following 
simple, but important result. 

\proclaim{Proposition 1.7} If $(X,S)$ is a symmetric set, then 
the map $J_n: X^n\to X^n$ 
given by the formula 
$$
J_n(x_1,...,x_n)=(f_{x_n}f_{x_{n-1}}...f_{x_2}(x_1),...,
f_{x_n}(x_{n-1}),x_n).\tag 1.16
$$
satisfies the commutation relation 
$$
J_n S^{ii+1}=\sigma^{ii+1} J_n.\tag 1.17
$$
\endproclaim

\demo{Proof} We will prove this statement by induction in $n$. 
For $n=2$, the statement follows directly from the involutivity of $S$. 
So let us assume the statement for $n=k$, and prove it for $n=k+1$. 

Observe that 
$$
J_{k+1}=Q_{k+1}\circ (J_k\times id_X),\tag 
1.18
$$
where $Q_n(x_1,...,x_n)=(f_{x_n}(x_1),...,f_{x_n}(x_{n-1}),x_n)$. 
Since $Q_n$ commutes with $\sigma^{ii+1}$ when $i<n-1$, 
formula (1.17) for $i<k$ follows from the induction 
assumption. So it remains to prove the formula for $i=k$. 

For $i=k$, the formula reduces to formulas (1.14),(1.15).
$\square$\enddemo

If in addition $(X,S)$ is nondegenerate, the map $J_n$ is obviously
bijective. Therefore, we get

\proclaim{Corollary 1.8} If $(X,S)$ is a nondegenerate symmetric set, 
then $J_n$ conjugates the twisted action of $S_n$ on $X^n$ 
to the canonical action of $S_n$ on $X^n$ by permutations. 
Thus, the two actions of $S_n$ are isomorphic. 
\endproclaim

Note that for a degenerate symmetric set the two actions of $S_n$ may be 
non-isomorphic. For example, for any set $X$ set $S(x,y)=(x,y)$.
Then $(X,S)$ is a symmetric set, but it is degenerate for $|X|>1$. 
As a result, the two actions of $S_n$ are not isomorphic in this case, since 
the twisted action of $S_n$ is trivial (all points are fixed). 

\head 2. The structure group\endhead 

\subhead 2.1. The structure group $G_X$ and its actions on $X$\endsubhead

Let $X$ be a set and $S:X^2\to X^2$ a mapping. It turns out to be 
very useful to introduce the group $G_X$ generated by elements 
of $X$ with defining relations
$$
xy=tz\text{ when }S(x,y)=(t,z).\tag 2.1
$$

\proclaim{Definition 2.1} The group $G_X$ is called the structure group of 
$X$. 
\endproclaim

{\bf Example.} If $(X,S)$ is the trivial pair 
($S(x,y)=(y,x)$), then $G_X=\Z^X$ is the free abelian group generated by $X$.  

One of 
the main properties of the structure group is the following: 

\proclaim{Proposition 2.1} Suppose that $(X,S)$ is nondegenerate. 
Then $(X,S)$ is a braided set
if and only if the following conditions are simultaneously satisfied:

(i)  the assignment $x\to f_x$ is a right action of 
$G_X$ on $X$;

(ii) the assignment $x\to g_x$ is a left action 
of $G_X$ on $X$;

(iii) the linking relation 
$$
f_{g_{f_y(x)}(z)}(g_x(y))=g_{f_{g_y(z)}(x)}(f_z(y))
$$
holds.  
\endproclaim

\demo{Proof} Conditions (i)-(iii) are exactly components 1-3
of the braid relation. 
$\square$\enddemo

\proclaim{Proposition 2.2} 
(a) Suppose $(X,S)$ is involutive, and the maps $f_x$ are invertible
and satisfy condition (i) of Proposition 2.1. Define the map $T:X\to
X$ by the formula $T(y)=f_y^{-1}(y)$. Then one has $f_x^{-1}T=Tg_x$.

Suppose in addition $g_x$ are invertible, so that
$(X,S)$ is nondegenerate and involutive. 
Then: 

(b) The map $T$ is invertible. Thus, the left actions of $G_X$ on $X$
given by $x\to f_x^{-1}$, $x\to g_x$ are isomorphic to each other.

(c) Condition (i) in Proposition 2.1 implies (ii) and (iii). 
Thus, $(X,S)$ is symmetric if and only if
the assignment $x\to f_x$ is a right action of $G_X$ on $X$. 
\endproclaim

\demo{Proof} 
(a) We have
$$
f_x^{-1}T(y)=f_x^{-1}f_y^{-1}(y)=
f^{-1}_{g_x(y)}f^{-1}_{f_y(x)}(y)=
f^{-1}_{g_x(y)}g_x(y)=Tg_x(y).
$$

(b) It follows from nondegeneracy and involutivity 
(see the 
proof of Proposition 1.4, move 1a)
that $T$ is invertible, and $T^{-1}(z)=g_z^{-1}(z)$. 

(c) The fact that (i) implies (ii) follows from (a) and (b). 
Let us prove that (i) implies (iii). Using (1.14), we can rewrite  
the linking relation in the form
$$
f_{f^{-1}_{f_zf_y(x)}(z)}f^{-1}_{f_y(x)}(y)=
f^{-1}_{f_{f_z(y)}f_{g_y(z)}(x)}f_z(y).\tag 2.2
$$
Set $u=f_zf_y(x)=f_{f_z(y)}f_{g_y(z)}(x)$
(these two expressions are equal by (i)). Then 
(2.2) takes the form
$$
f_{f^{-1}_u(z)}f^{-1}_{f_z^{-1}(u)}=
f^{-1}_uf_z,
$$
which is a direct consequence of (1.14) and (i). 
$\square$\enddemo

\subhead 2.2 The properties of the group $G_X$\endsubhead

Now we will determine the structure 
of the group $G_X$ for a nondegenerate symmetric set. 

Let $Aut(X)$ be the group of permutations of $X$, and $\Z^X$ be the 
free abelian group spanned by $X$. We will denote the generator
of $\Z^X$ corresponding to $x\in X$ by $t_x$. 
Let $M_X=Aut(X)\ltimes \Z^X$
be the semidirect product, associated to the action of $Aut(X)$ on $\Z^X$. 
The group $M_X$ consists of elements of the form 
$s t$, where $s\in Aut(X)$, $t\in \Bbb Z^X$, and we have the commutation relation 
$st_x=t_{s(x)}s$.

Consider the assignment
$$
x\to f_x^{-1}t_x.\tag 2.3
$$

\proclaim{Proposition 2.3} If $(X,S)$ is a nondegenerate 
symmetric set then assignment (2.3) extends to a group homomorphism 
$G_X\to M_X$. 
\endproclaim

\demo{Proof} We have to show that 
$$
f_x^{-1}t_xf_y^{-1}t_y=f_u^{-1}t_uf_v^{-1}t_v\text{ when } 
S(x,y)=(u,v).\tag 2.4
$$
We have
$$
f_x^{-1}t_xf_y^{-1}t_y=f_x^{-1}f_y^{-1}t_{f_y(x)}t_y,\
f_u^{-1}t_uf_v^{-1}t_v=f_u^{-1}f_v^{-1}t_{f_v(u)}t_v.
$$
But we have $v=f_y(x)$, and by involutivity of $S$ we have $f_v(u)=y$.
Therefore, (2.4) follows from (1.15). 
$\square$\enddemo

We will denote the constructed homomorphism $G_X\to M_X$ by $\phi_f$. 
Thus, $\phi_f(x)=f_x^{-1}t_x$. 

Analogously, we can construct a homomorphism $\phi_g:G_X\to M_X$ 
given by $\phi_g(x)=t_xg_x$. The homomorphisms $\phi_f,\phi_g$ 
are conjugate in the following sense: if we denote by 
$\hat T$ the automorphism of $M_X$ induced by the permutation 
$T$ of $X$, then $\hat T$ conjugates $\phi_f$ to $\phi_g$. 
Thus, it is enough for us to study the properties of $\phi_f$.  
 From now on we will denote it simply by $\phi$. 
\proclaim{Proposition 2.4} 
 The homomorphism $\phi$ is injective. 
\endproclaim
\demo{Proof}
Proposition 2.4 is an immediate corollary of Proposition 2.5 (see below).
$\square$\enddemo
Let $\pi:G_X\to \Z^X$ be the map
defined by $\pi(g)=t$ if $\phi(g)=st$, $s\in Aut(X)$, $t\in \Z^X$.
We fix the structure of a left $G_X$-module on $\Z^X$ induced from the
assignment
$x\to f_x^{-1}$.

For convenience, we will write the group operation in $\Z^X$ additively. 

\proclaim{Proposition 2.5}
(a) $\pi$ is a 1-cocycle of $G_X$ with coefficients in the $G_X$ -  module
$\Z^X$,
i.e. $\pi(g_1g_2)=g_2^{-1}\pi(g_1)+\pi(g_2)$. 

(b) $\pi$ is bijective.
\endproclaim
The proof of Proposition 2.5 is contained in the next section. 

\subhead 2.3. Proof of Proposition 2.5\endsubhead
Property (a)  follows from the definition of the semidirect
product and the fact that $\phi$ is a homomorphism. So, we have to prove
(b). We will explicitly construct the map
$h:\Z^X\to G_X$ inverse to $\pi$. Let 
$X^+=\{t_x\in \Z^X| x\in X\}$, $X^-
= \{-t_x\in \Z^X|
x\in X \}$
and $Y=X^+\cup X^- $. They are  $G_X$-invariant subsets of $\Z^X$.
For a nonnegative integer k consider the subset $\Z_k^X$ of $\Z^X$
consisting
of all the elements that are sums of no more than k elements of $Y$.
In particular, $\Z_0^X=\{0\}$, $\Z_k^X\subset \Z_{k+1}^X$ and $\Z_k^X$ form
a covering of $\Z^X$. Similarly, let $G_X^k$ denote the set of elements
of $G_X$ representable as a product  $x_1...x_m$, where 
$x_1,...,x_m\in X\cup X^{-1}\subset G_X$ and $m\leq k$. 
%It is clear from
%(a) that $\pi$ maps  $G_X^k$ to $Z_k^X$. 

 We want to define $h$ inductively on each $\Z_k^X$
in a compatible way. Let $h(0)=1$, $h(t_x)=x,\ h(-t_x)=(g_x^{-1}(x))^{-1}$
for $x\in X$. In this way, $h$ is defined on $\Z_1^X=Y$.

For convenience, for $g\in G_X$ and $\xi\in \Z^X$, denote by $g*\xi$ 
the result of the action of $g$ on $\xi$. 

\proclaim{Lemma 2.6}
For $\xi,\eta \in Y$ one has $ h(h(\xi)*\eta)h(\xi)= h(h(\eta)*\xi)h(\eta)$.
\endproclaim
\demo{Proof}

We have to consider 3 cases: 

(i) Both $\xi$ and $\eta$ belong to $X^+$.
Let $x=h(\xi),y=h(\eta)$. 
We want to show that
$f_x^{-1}(y)x=f_y^{-1}(x)y$. Let $S(f_x^{-1}(y),x)=(z, y)$, then
$S(z, y)=(f_x^{-1}(y),x)$, so $f_y(z)=x$ and hence
$S(f_x^{-1}(y),x)=(f_y^{-1}(x),y)$.  

(ii) Only $\xi$ belongs to $X^+$, while 
$\eta\in X^-$. Then we have $h(\xi)=x\in X$, and 
$\eta=-t_{y'}$ for some $y'\in X$.  
We need to check that 
$(g_{f_x^{-1}(y')}^{-1}(f_x^{-1}(y')))^{-1}x=
f_{g_{y'}^{-1}(y')}(x)(g_{y'}^{-1}(y'))^{-1}$. Let
$z=g_{y'}^{-1}(y')=T^{-1}(y')$ (see the proof of Proposition 2.2 (a),
(b)). Then $g_{f_x^{-1}(y')}^{-1}(f_x^{-1}(y'))=
T^{-1}f_x^{-1}T(z)=g_x(z)$. So, the desired equality is just 
$(g_x(z))^{-1}x=f_z(x)z^{-1}$, which holds by the definition of $G_X$.

(iii) Both $\xi$ and $\eta$ belong to $X^-$. This case is similar to
(ii), so the proof is omitted.  
$\square$\enddemo

Now, let us assume that $h$ has already been defined for elements of
$\Z_k^X$. Take $\eta\in \Z_{k+1}^X$, then $\eta=a+\xi$ for $a\in
\Z_k^X$, $\xi\in Y$. Define $h(\eta)=h(h(a)*\xi)h(a)$.

\proclaim{Lemma 2.7} The map $h$ is well-defined on each $\Z_k^X$, and thus
on the
whole $\Z^X$. 
\endproclaim

\demo{Proof} 
We proceed by induction on $k$.
The lemma is certainly true for $k=0,1$. Suppose $h$ is well-defined on
$\Z_{k-1}^X$.
For any $a\in \Z_{k-2}^X$ and $\xi,\eta \in Y$ 
$$h((a+\xi)+\eta)=h(h(a+\xi)*\eta)h(a+\xi)=
h(h(h(a)*\xi)h(a)*\eta)h(h(a)*\xi)h(a).$$ 
We need to check that the last expression is symmetric in $\xi, \eta$ and
equal to $h(a)$ when $\xi=-\eta$. Set $\xi'=h(a)*\xi,\
\eta'=h(a)*\eta$, then we arrive at the formula \ \ \ 
$h((a+\xi)+\eta)=h(h(\xi')*\eta')h(\xi')h(a)$ that is symmetric 
in $\xi', \eta'$ (hence in $\xi, \eta$) 
by Lemma 2.6. If $\xi'=-\eta'=t_z$ for $z\in X$
we have
$$h(z*(-t_z))zh(a)=h(-t_{T(z)})zh(a)=(T^{-1}T(z))^{-1}zh(a)=h(a).$$ 
Lemma 2.7 is proved.
$\square$\enddemo
%It follows from the definition of $h$ that it maps every $\Z_k^X$ to
%$G_X^k$.   

\proclaim{Lemma 2.8} The maps $h:\Z^X\to G_X$ and $\pi:G_X\to \Z^X$
are inverse to each other. 
\endproclaim
\demo{Proof}
It suffices to check that:

(i)$\pi\circ h(a)=a$ for $a\in \Z_k^X$.

(ii)$h\circ \pi(b)=b$ for $b\in G_X^k$.

The statements (i), (ii) are quite simple for $k=0,1$. Let us make an
inductive step. For $a\in \Z_k^X,\ \xi\in Y$ 
$$\pi\circ
h(a+\xi)=\pi(h(h(a)*\xi)h(a))=h(a)^{-1}*\pi(h(h(a)*\xi)) + 
\pi(h(a))=a+\xi,$$ provided $\pi\circ h(a)=a$. So, (i) is proved. Similarly,
for $b\in G_X^k,\ y\in X\cup X^{-1}$
$$h \circ\pi(yb)=h(b^{-1}*\pi(y)+\pi(b))=h(b^{-1}*(b*\pi(y)))b=yb,$$ 
provided $h \circ\pi(b)=b$. So, Lemma 2.8 is proved along with Proposition
2.5.
$\square$\enddemo
\subhead 2.4. Classification of nondegenerate symmetric sets 
via groups with bijective 1-cocycles
\endsubhead

\proclaim{Definition 2.2} A bijective cocycle quadruple 
is a quadruple $(G,X,\rho,\pi)$, where $G$ is a group, $X$ is a set, 
$\rho:G\times X\to X$ a left action of $G$ on $X$, and 
$\pi:G\to \Z^X$ a {\it bijective} 1-cocycle of $G$ with coefficients in 
$\Z^X$, where $G$ acts in $\Z^X$ by $\rho$. 
\endproclaim

\proclaim{Theorem 2.9} Nondegenerate symmetric sets, up to 
isomorphism, are in 1-1 correspondence with bijective cocycle 
quadruples, up to isomorphism. This correspondence 
is given by $F(X)=(G_X,X,\rho,\pi)$, where $\rho$ is the action 
of $G_X$ on $X$ by $x\to f_x^{-1}$, and $\pi$ is the 1-cocycle of 
Proposition 2.5. 
\endproclaim

In fact,  nondegenerate symmetric sets form a category,
as well as bijective cocycle 
quadruples. Morphisms in these categories are just maps 
which preserve all structures. Namely, a morphism of nondegenerate 
symmetric sets $(X_1,S_1)\to (X_2,S_2)$ is a map 
$f:X_1\to X_2$ such that $(f\times f)S_1=S_2(f\times f)$. 
Similarly, a morphism of bijective cocycle quadruples
$(G_1,X_1,\rho_1,\pi_1)\to (G_2,X_2,\rho_2,\pi_2)$ is a pair of maps 
$f:X_1\to X_2,\phi:G_1\to G_2$, where $\phi$ is a group homomorphism, 
and $\rho_2(\phi(g))\circ f=f\circ \rho_1(g)$, 
$\pi_2(\phi(g))=f(\pi_1(g))$ for any $g\in G_1$
(in the last formula $f$ is extended to a group homomorphism
$\Z^{X_1}\to \Z^{X_2}$). 
Denote these categories by 
NSS and BCQ (by abbreviating the names). It is clear 
that the map $F$ is not only a map but also a functor 
$NSS\to BCQ$. Theorem 2.9 can be strengthened as follows.

\proclaim{Theorem 2.10} The functor $F:NSS\to BCQ$ is an equivalence 
of categories. 
\endproclaim

\demo{Proof of Theorems 2.9,2.10} To prove Theorems 2.9, 2.10   
 it is necessary to construct the inverse functor to the functor $F$, 
i.e. learn to reconstruct $(X,S)$ from the quadruple $F(X,S)$.

Consider a bijective cocycle quadruple $(G,X,\rho,\pi)$.
Let $t_x$ be the generator of $\Z^X$ 
corresponding to $x\in X$.
Then we have a natural embedding $a:X\to G$ given by the formula
$a_x=\pi^{-1}(t_x)$. 

For any $x\in X$, define the map $f_x:X\to X$ by 
$f_x^{-1}=\rho(a_x)$. 
Define the map $g_x:X\to X$ by $g_x(y)=f^{-1}_{f_y(x)}(y)$. 

Let us show that 
$$
a_xa_y=a_{g_x(y)}a_{f_y(x)}.\tag 2.5
$$ 
Indeed, 
we have
$$
\pi(a_xa_y)=\rho(a_y)^{-1}\pi(a_x)+\pi(a_y)=\rho(a_y)^{-1}t_x+t_y=
t_{f_y(x)}+t_y.
$$
Similarly, 
$$
\pi(a_{g_x(y)}a_{f_y(x)})=t_{f_{f_y(x)}g_x(y)}+t_{f_y(x)}=
t_y+t_{f_y(x)}.
$$
Since $\pi$ is bijective, we have the desired equality 
$a_xa_y=a_{g_x(y)}a_{f_y(x)}$.

Let us show that $g_x$ is invertible. For this purpose define 
the map $T:X\to X$ by the formula $T(x)=f_x^{-1}(x)$. It is enough 
for us to show that $T$ is invertible, since by (2.5) 
(similarly to Proposition 2.2(a))
we have $f_x^{-1}T(y)=Tg_x(y)$. 

Proof that $T$ is invertible. We have $t_{T(x)}=
t_{f_x^{-1}(x)}=\rho(a_x)(t_x)=\rho(a_x)(\pi(a_x))$. Since 
$\pi$ is a 1-cocycle, $\rho(a_x)(\pi(a_x))=-\pi(a_x^{-1})$.
It is clear that the inverse of T is given by the map 
$T':X\to X$, $T'(x)=(\pi^{-1}(-t_x))^{-1}$ if the latter is well defined,
i.e. we have to check that $T'(x)\in X$ for $x\in X$. One observes that
by the definition of $\pi$ 
$\pi(\rho(\pi^{-1}(-t_x))(x)\pi^{-1}(-t_x))=0$ and hence
$(\pi^{-1}(-t_x))^{-1}= \rho(\pi^{-1}(-t_x))(x)\in X$. So, $T$ is
bijective.

Now define a map $S:X^2\to X^2$ by the formula 
$S(x,y)=(g_x(y),f_y(x))$. By the construction, it is 
nondegenerate and involutive. 
The fact that the map $S$ 
satisfies condition (i) of Proposition 2.1 follows from (2.5). 

Thus, $(X,S)$ is a nondegenerate symmetric set. In other words, 
we constructed a map $F':BCQ\to NSS$, by 
$F'(G,X,\rho,\pi)=(X,S)$. 

It is clear that the map $F'$ is a functor. 
Indeed, if $\mu:(G,X,\rho,\pi)\to (G',X',\rho',\pi')$ is a morphism 
of bijective cocycle quadruples, then it respects the assignment 
$x\to a_x$, and therefore respects $f$ and $g$, so it defines a morphism 
$F'(\mu):F'(G,X,\rho,\pi)\to F'(G',X',\rho',\pi')$.  

To complete the proof, we need to show that 
$F\circ F'=id,F'\circ F=id$. 

The identity $F'\circ F=id$ is obvious. Let us prove that
$F\circ F'=1$. 

Let $(G,X,\rho,\pi)$ be a bijective cocycle quadruple, and 
$(X,S)=F'(G,X,\rho,\pi)$. We need to show that $G=G_X$, 
and $\rho,\pi$ are defined in the standard way. 

As we have seen, 
the map $X\to G$ defined by $x\to a_x$ extends to a homomorphism 
$a: G_X\to G$, which transforms $\rho$ into the standard 
action of $G_X$ on $X$ ($x\to f_x^{-1}$)
and the composition $\pi\circ a$ is the standard 
bijective 1-cocycle on $G_X$. Since $\pi$ is bijective, we get that $a$
is bijective, as desired. 

Theorems 2.9, 2.10 are proved. 
$\square$\enddemo

{\bf Remark 1.} Groups with bijective 1-cocycles have a nice geometric 
interpretation, which was pointed out to us by R.Howe and G.Margulis.
Namely, in the case of Lie groups, a bijective 1-cocycle on a simply connected 
group $G$ with coefficients in some real representation is the same 
thing as a left-invariant affine structure on $G$ as a manifold, which 
identifies $G$ with an affine space. Of course, to admit such a structure, 
the group must be contractible topologically, i.e. solvable. We will show 
later that this is also the case for finite groups with a bijective 
1-cocycles. 
  
{\bf Remark 2.} It is known that many solvable Lie groups admit left-invariant 
affine structures as above \cite{Bu}. It was even believed that any solvable 
simply connected Lie group does (Milnor's conjecture), but it was 
recently disproved by D.Burde \cite{Bu}, who found a nilpotent group 
which has no such structures. 

{\bf Remark 3.} Bijective cocycles on finite groups were recently used 
by the first author and S.Gelaki to construct new examples of semisimple
Hopf algebras \cite{EG}. 

\subhead 2.5. The canonical abelian subgroup $\Gamma$ of $G_X$\endsubhead 
 
Let $\Gamma\subset G_X$ be the intersection $G_X\cap \Z^X$, where 
$G_X$ and $\Z^X$ are both regarded as subgroups of $M_X$.

It is clear that $\Gamma$ is a normal subgroup of $G_X$. 
Indeed, $\Gamma$ is the kernel of the homomorphism 
$\rho:G_X\to Aut(X)$, defined by the action of $G_X$ on $X$. 

Let $G_X^0=G_X/\Gamma$, and $A=\Z^X/\Gamma$. 
The group $G_X^0$ is the image of $G_X$ in $Aut(X)$,  so, in particular, 
it is finite if the set $X$ is finite. 

The 1-cocycle $\pi$ is obviously equal to $id$ on $\Gamma$, so it descends to
a bijective 1-cocycle $\bar\pi: G_X^0\to A$. This follows from the fact that
$\pi(b\gamma)=\pi(b)+\gamma$, $b\in G_X$, $\gamma\in \Gamma$. 

If $X$ is finite, $A$ is a finite abelian group, and $|A|=|G_X^0|$.  The resulting quadruple $(G_X^0,A,\bar\rho,\bar\pi)$ where $\bar\rho: G_X^0 \rightarrow Aut(A)$ descends from $\rho$ motivates the following general definition:

\proclaim{Definition 2.3} Define a bijective cocycle datum to be a quadruple
$(G,A,\rho,\pi)$ where $G$ is a group, $A$ a $G$-module with the $G$-action 
given by the homomorphism $\rho: G \rightarrow Aut(A)$, and 
$\pi:G \rightarrow A$ a bijective 1-cocycle employing the action $\rho$.  
\endproclaim

In fact, given any bijective cocycle datum, if we choose any set $X$
with a $G$-action $\rho': G \rightarrow Aut(X)$ and a map $\phi: X
\rightarrow A$ carrying $\rho'$ to $\rho$, then $X$ carries the
structure of a nondegenerate set. Explicitly, define
$f_y(x)=\rho'(\pi^{-1}(y)^{-1})(x)$, $g_x(y)=f^{-1}_{f_y(x)}(y)$, and
$S(x,y)=(g_x(y),f_y(x))$ for $x,y \in X$. Then $(X,S)$ is a
nondegenerate symmetric set.  Let us therefore make the following
definition.

\proclaim{Definition 2.4} Take any bijective cocycle datum $(G,A,\rho,\pi)$. Define
a set-structure for the datum to be a triple $(X,\rho',\phi)$ where
$X, \rho', \phi$ are as above.
\endproclaim

Note that in any set-structure, if $\phi(X)$ does not generate $G$,
then we could consider the subgroup $G' \subset G$ generated by $X$,
and clearly $\pi(G')$ is generated by $\pi(X)$ yielding the smaller
datum $(G',A',\rho,\pi)$.  Thus, let us call a set-structure {\it
generating} if $\phi(X)$ generates $G$; these are the kinds we will
consider.  Note that when $\phi(X)$ generates $G$ freely the datum is
merely the corresponding bijective cocycle quadruple to the induced
solution $(X,S)$.

Now we may generalize the notion of factoring by $\Gamma$ as follows:
Take any datum $(G,A,\rho,\pi)$ together with a generating
set-structure $(X,\rho',\phi)$.  Let $H=Ker(\rho')$, $K=\pi(H)$.
Since $X$ is generating, $H$ is a subgroup of $Ker(\rho)$ so $K$ is a
subgroup of $A$ since $\pi$ is an injective homomorphism restricted to
$H$.  Then, form $G/H, A/K$.  Now $G/H$ is just the image of $G$ in
$Aut(X)$.  Then, let $\bar\rho:G/H \rightarrow Aut(A/K)$, $\bar\pi:
G/H \rightarrow A/K$ be the maps descending from $\rho, \pi$. We
arrive at the conclusion that $(G/H,A/K,\bar\rho,\bar\pi)$ is a
bijective cocycle datum.  Together with the same set structure (taking
the images of $\rho'$ and $\phi$ in $G/H$) we obtain the solution
$(X,S)$ which is isomorphic to the original solution since the actions
$f_x^{-1}:X \rightarrow X, x \in X$ are the same.  It is clear that
this is the smallest such datum giving rise to this solution.

Call a generating set-structure {\it faithful} if $G$ acts faithfully
on $X$. Then we have the following result:

\proclaim{Proposition}  Given a bijective cocycle datum $(G,A,\rho,\pi)$ and 
a faithful generating set-structure $(X,\rho',\phi)$, then 
$(G,A,\rho,\pi)$ is isomorphic to $(G_X^0,A_X,\rho_X,\pi_X)$ where $G_X^0, A_X, \rho_X$ and
$\pi_X$ are obtained from the solution $(X,S)$ induced by the original datum
and set-structure.
\endproclaim

\demo{Proof}  Indeed, it is clear by the construction in Section 2.4 that
$xy=g_x(y) f_y(x)$ in $G$ so $G$ is the image of $G_X$ under the homomorphism
sending $x \in X \subset G_X$ to $\phi(x)$ in $G$.  Since $G$ acts faithfully
on $X$, it must be isomorphic to the image of $G_X$ in $Aut(X)$, namely $G_X^0$.
The action $\rho$ is determined by $\rho'$ and hence is the same as $\rho_X$ since
$\phi(X)$ generates $G$, so it is easy to see that $A$ and the cocycle $\pi$ are
the same as $A_X$ and $\pi_X$.
\enddemo

There are two main ideas behind this development.  First, we have seen
that the bijective cocycle structure can be formulated without the
choice of a set $X$, and this leads easily to nondegenerate symmetric
sets once $X$ is chosen.  Secondly, we now may classify nondegenerate
symmetric sets via faithful generating set-structures on such data,
which is especially useful in the case where $X$ is finite and hence
so is $G_X^0$.  The usefulness of this formulation will be firmly
established in the following section.

\subhead 2.6. Indecomposable symmetric sets of prime order
\endsubhead

\proclaim{Definition 2.5} (a) A subset $Y$ of a nondegenerate 
symmetric set $X$ is said to be an invariant subset if
$S(Y\times Y)\subset Y\times Y$.

(b) An invariant subset $Y\subset X$ is said to be nondegenerate 
if $(Y,S|_{Y\times Y})$ is a nondegenerate symmetric set.
 
(c) A nondegenerate symmetric set $(X,S)$ 
is said to be decomposable if it is a union of two nonempty disjoint
nondegenerate invariant subsets. Otherwise, $(X,S)$ is said to be 
indecomposable.
\endproclaim

For example, a permutation solution is indecomposable if and only 
if it is cyclic.

{\bf Remark.} If $X$ is finite, then any invariant subset $Y$ of $X$ is 
nondegenerate. Indeed, the map $S|_{Y\times Y}$ is bijective, as $S^2=1$, 
and for any $y\in Y$, the maps 
$f_y,g_y:Y\to Y$ are injective (by nondegeneracy of $X$), hence bijective
by the finiteness. 
Thus $Y$ is nondegenerate. 

\proclaim{Proposition 2.11} A nondegenerate symmetric set $(X,S)$ is 
indecomposable if and only if $G_X$ acts transitively on $X$. 
\endproclaim

\demo{Proof} Proposition 2.15 below implies 
 that if $G_X$ acts transitively on $X$ then $(X,S)$
is indecomposable. Conversely, if $x \to f_x^{-1}$ is  not
transitive, consider two complementary nonempty $G_X$-invariant subsets
$X_1$ and $X_2$.  They are invariant under $f_x$ for all $x$, and hence under
$T$ (since $T(y)=f_y^{-1}(y)$), so they are invariant under $g_x$ for all
$x$.  Thus $X_1$ and $X_2$ are invariant subsets of $X$. 
It is clear that these subsets are nondegenerate. Thus $(X,S)$ is decomposable.
$\square$\enddemo

Now we will classify finite indecomposable nondegenerate symmetric sets 
which have $p$ elements, where $p$ is a prime. 

\proclaim{Theorem 2.12} Let $(X,S)$ be an indecomposable nondegenerate 
symmetric set, and $|X|=p$, where $p$ is a prime. Then 
$(X,S)$ is isomorphic to the cyclic permutation solution 
$(\Z/p\Z,S_0)$, where $S_0(x,y)=(y-1,x+1)$. 
\endproclaim

\demo{Proof} Since by Proposition 2.11
the group $G_X^0$ acts transitively on $X$, its order 
is divisible by $p$. Since $G_X^0\subset Aut(X)$, its order divides $p!$. 
Thus, $|G_X^0|=|A|=pn$, where $n$ is coprime to $p$. 

Thus, $A=\Z/p\Z\oplus A_0$, where $|A_0|$ is coprime to $p$. 

The subgroup $A_0\subset A$ is the group of all elements of order 
not divisible by $p$. Therefore, $A_0$ is $G_X^0$-stable, and 
the cocycle $\bar\pi:G_X^0\to A$ defines a 1-cocycle 
$\pi':G_X^0\to \Z/p\Z$. 

Set $H=(\pi')^{-1}(0)=\bar\pi^{-1}(A_0)$. It is easy to see that $H$ is a subgroup 
of $G_X^0$ of order $n$. 

Let $H_x$ be the stabilizer of a point $x\in X$ in $G_X^0$. This is a subgroup 
of $G_X^0$ of index $p$, i.e. of order $n$. 
We want to show that $H_x=H$ for all $x$. 
Since $|H_x|=|H|$, 
for this purpose it is enough to show that $\pi'(H_x)=\{0\}$. 

We will show that $H_x$ acts trivially on $A/A_0=\Z/p\Z$. 
This implies that $\pi'|_{H_x}$ is simply a homomorphism $H_x\to \Z/p\Z$.
This, in turn, implies that $\pi'|_{H_x}=0$, since $|H_x|=n$ is 
coprime to $p$. 

To show that $H_x$ acts trivially on $A/A_0$, it is enough to prove that 
the image $\bar t_x\in A/A_0$ of the element $t_x\in \Z^X$ is not zero. 
Indeed, in this case $\bar t_x$ generates $A/A_0$, while $\bar t_x$ is 
by definition fixed by $H_x$. 

Now we prove that $\bar t_x\ne 0$. Assume that $\bar t_x=0$. 
Since $G_X$ acts transitively on $X$, we get $\bar t_y=0$ for all $y\in X$. 
Thus, the natural map $\Z^X\to A/A_0$ is zero. Contradiction. 

Thus, we showed that $H_x=H$. Therefore, $H$ acts trivially on $X$. This 
implies that $H=A_0=\{1\}$, $G_X^0=A=\Z/p\Z$, the action of $G_X^0$ on $A$ 
is  trivial, and $\bar\pi=id$. 

To conclude the proof, it is enough to observe that 
$\bar t_x=\bar t_y$ for any $x,y\in X$. This follows from the fact 
that there exists $g\in G_X^0$ such that $gx=y$, while the action 
of $G_X^0$ on $A$ is trivial. Thus, $\bar t_x$ does not depend on $x$, and 
hence $a_x=\bar\pi^{-1}(\bar t_x)\in G_X^0$ does not depend on $X$. 
Therefore, $(X,S)$ is a permutation solution. The theorem is proved.
$\square$\enddemo
\subhead 2.7. Solvability of the structure group
\endsubhead

In this section we prove the solvability of the structure group
of a finite nondegenerate symmetric set.

Let $G$ be a finite group, $p$ a prime divisor 
of $|G|$. Let us write $|G|=mp^k$
for a  positive $k$ and $m$ coprime to $p$.
\proclaim{Definition 2.6}
A subgroup $H$ of $G$  of order m is called a Hall $p'$-subgroup.
\endproclaim
  
We will use the following theorem of Hall: 

\proclaim{Theorem 2.13}\cite{As}
If a finite group $G$ has a Hall $p'$-subgroup for each prime $p$
dividing $|G|$, then it is solvable.
\endproclaim
\proclaim{Theorem 2.14}
The structure group $G_X$ of a finite nondegenerate symmetric set is
solvable. 
\endproclaim
\demo{Proof}
It is enough to show that the finite group $G^0_X=G_X/\Gamma$ is solvable.
We shall make use of the 1-cocycle $\bar\pi:G^0_X\to A$ defined in section
2.5. Let $|G^0_X|=mp^k$, $m, p$ being coprime. Define
$H_p=\bar\pi^{-1}(p^kA)$. Obviously, $p^kA$ is invariant under the action
of 
$G^0_X$, so $H_p$ is a subgroup of $G^0_X$. The order of $H_p$
equals the order of $p^kA$, that is, in turn, equals $m$. So, $H_p$
is a Hall $p'$-subgroup of $G^0_X$ and we apply Lemma 2.13 to conclude
that
$G^0_X$ is solvable.  
$\square$\enddemo

\subhead 2.8. Extensions and the structure group \endsubhead

Let $(Z,S)$ be a nondegenerate symmetric set. 

\proclaim{Definition 2.7} We will say that 
$Z$ is a union of two nondegenerate 
invariant subsets $X,Y$ if $X\cap Y=\emptyset$ and
$Z=X\cup Y$ as a set. 
\endproclaim

Thus, $Z$ can be represented as a union of two nonempty 
nondegenerate invariant subsets if and only if it is decomposable. 

\proclaim{Proposition 2.15} If $(Z,S)$ is 
a union of nondegenerate 
invariant subsets $X$ and $Y$, then the map $S$ defines 
bijections $X\times Y\to Y\times X$, and 
$Y\times X\to X\times Y$. 
\endproclaim

\demo{Proof} It is clear that $S$ defines an automorphism 
of $X\times Y\cup Y\times X$. So we need to show that it does not map 
elements of $X\times Y$ to $X\times Y$, and the same for $Y\times X$. 

Let $x\in X,y\in Y$, and assume that $S(x,y)=(x',y')$, $x'\in X,y'\in Y$. 
Thus, $g_x(y)=x'$. However, by nondegeneracy of $X$ there exists 
$x''\in X$ such that $g_x(x'')=x'$. This violates the nondegeneracy of $Z$.
Contradiction. $\square$\enddemo

{\bf Example.} The most obvious example of a union is 
the trivial union, defined by $S(x,y)=(y,x)$, and 
$S(y,x)=(x,y)$. However, as we will
see later, there are much more interesting ways of constructing unions. 

Let $(X,S_X),(Y,S_Y)$ be nondegenerate symmetric sets. 
Denote by $Ext(X,Y)$ (extensions of $X$ by $Y$)
the set of all decomposable solutions 
$Z$ which are unions of $X$ and $Y$. As we showed,
an element $Z\in Ext(X,Y)$ is completely determined by the function 
$S_Z:X\times Y\to Y\times X$. 

Let us write $S_Z(x,y)$ in the form $(g_x(y),f_y(x))$. 

\proclaim{Proposition 2.16} If $Z\in Ext(X,Y)$ then 
the assignments $x\to g_x$, $y\to f_y^{-1}$ are actions of 
$G_X$ on $Y$ and of $G_Y$ on $X$. 
\endproclaim

\demo{Proof} The statement follows from Proposition 2.1, 
since $G_X,G_Y$ are obviously subgroups of $G_Z$. 
$\square$\enddemo

The conclusion of Proposition 2.16 is clearly not sufficient 
for $S_Z$ to define an extension. There is an additional, rather complicated 
linking condition between $f_y$ and $g_x$. Thus, we will consider 
a special case of ``one-sided extensions''.

\proclaim{Definition 2.8} An element $Z\in Ext(X,Y)$ is called 
a right (respectively, left) extension of $X$ by $Y$ if 
$S_Z(x,y)=(y,f_y(x))$ (respectively, 
$S_Z(x,y)=(g_x(y),x)$) for $x\in X,y\in Y$. 
The set of right (left) extensions of $X$ by $Y$ will be denoted by
$Ext_+(X,Y)$, $Ext_-(X,Y)$, respectively. 
\endproclaim

It is clear that $Ext_+(X,Y)=Ext_-(Y,X)$. 

The following proposition gives a complete group-theoretic description of
$Ext_+(X,Y)$. 

\proclaim{Proposition 2.17} The formula 
$S_Z(x,y)=(y,f_y(x))$ defines an element $Z\in Ext_+(X,Y)$ 
if and only if the assignment $y\to f_y^{-1}$ defines an action 
of $G_Y$ on $(X,S_X)$. 
\endproclaim

\demo{Proof}
We have
$$
S^{12}S^{23}S^{12}(x_1,x_2,y)=
(y,f_yg_{x_1}(x_2),
f_yf_{x_2}(x_1)),\tag 2.6
$$
and 
$$
S^{23}S^{12}S^{23}(x_1,x_2,y)=
(y,g_{f_y(x_1)}f_y(x_2),
f_{f_y(x_2)}f_y(x_1)).\tag 2.7
$$
Equating (2.6) and (2.7) 
shows that $f_y$ preserves $S_X$. It is easy to check 
that other relations for $Z$ do not impose any new restrictions. 
The proposition is proved. 
$\square$\enddemo

\proclaim{Corollary 2.18} If $Z\in Ext_+(X,Y)$ then the group $G_Z$ 
is isomorphic to $G_Y\ltimes G_X$, where the semidirect product 
is formed using the action of $G_Y$ on $X$ via $y\to f_y^{-1}$. 
\endproclaim

In particular, for the trivial union $Z=X\cup Y$ 
we have $G_Z=G_X\times G_Y$. 

Finally, let us describe the set $Ext(X,Y)$, where 
$X,Y$ are trivial symmetric sets, i.e. $S_X,S_Y$
are the permutation of components. In this case, 
$G_X=\Z^X$, and $G_Y=\Z^Y$. 

\proclaim{Proposition 2.19} The formula 
$S_Z(x,y)=(g_x(y),f_y(x))$ defines an element 
$Z\in Ext(X,Y)$ if and only if 

(i) the assignments $x\to g_x$, $y\to f_y^{-1}$ are actions of 
$\Z^X$ on $Y$ and of $\Z^Y$ on $X$;

(ii) the homomorphisms $\rho_X:\Z^X\to Aut(Y)$, $\rho_Y:\Z^Y\to Aut(X)$
defined by (i) are invariant under $\Z^Y$, $\Z^X$, respectively
(where the actions of $\Z^X$ on $Aut(X)$ and of $\Z^Y$ on $Aut(Y)$ 
are trivial, and the actions of $\Z^X,\Z^Y$ on each other are by $\rho_X$, 
$\rho_Y$).  
\endproclaim

\demo{Proof} Straightforward.
$\square$\enddemo

\subhead 2.9. The 
quantum algebras associated to a nondegenerate symmetric set, and their 
relation to the structure group  
\endsubhead

Let $(X,S)$ be a finite nondegenerate symmetric set.   
Let $V$ be the complex vector space spanned by $X$.
Let $v_x$ be the vector in $V$ corresponding to $x\in X$, and 
$E_{xy}:V\to V$ be the endomorphism of $V$ defined by 
$E_{xy}v_z=\delta_{yz}v_x$. 

Let $R=\sigma S$. We regard
$R$ as a linear operator $V\o V\to V\o V$. 

Following Faddeev, Reshetikhin, Sklyanin, and Takhtajan \cite{FRT},
define two quadratic algebras over $\C$ associated to 
$X$.

1. The quantized algebra of functions on $V$.  
This is the quadratic algebra $Q_X$ with generators
$q_x$, $x\in X$, and relations 
$$
R^{12}q^{13}q^{23}=q^{23}q^{13},
\tag 2.8
$$
where 
$$
q:=\sum_x v_x\o q_x\in V\o Q_X\tag 2.9
$$ 

2. The quantized algebra of functions on $\text{End}(V)$. 
This is the quadratic algebra $A_X$ with generators 
$L_{xy}$, $x,y\in X$, and relations
$$
R^{12}L^{13}L^{23}=L^{23}L^{13}R^{12},\tag 2.10
$$ 
where $L=\sum E_{xy}\o L_{xy}$, $E_{xy}\in End(V)$.

{\bf Remark.} The algebra $A_X$ is a special case of $H_R$ of \cite{FRT};
similar algebras were also studied in \cite{Ha,Sch}. 

Observe that the relations of $Q_X$ can be written in the form
$$
q_xq_y=q_{g_x(y)}q_{f_y(x)}.\tag 2.11
$$

Let $G_X^+$ be the set of elements of $G_X$ representable as a product 
of the generators (without inverses). This set is a monoid.  
We have $G_X^+=\cup_{n\ge 0}G_X^{+n}$, where $G_X^{+n}$ is the set of elements 
representable as a product of $n$ generators. It follows from 
Proposition 2.5 that this is a $\Z_+$-grading of $G_X^+$, 
and $|G_X^{+n}|=\left(\matrix
n+N-1\\ n\endmatrix\right)$. 

Equation (2.11) implies that the algebra $Q_X$ is isomorphic 
to $\C[G_X^+]$. 

Similarly, the relations in $A_X$ can be written in the form
$$
L_{xz}L_{yt}=L_{g_x(y)g_z(t)}L_{f_y(x)f_t(z)}.\tag 2.12
$$
Thus, $A_X$ is isomorphic to $\C[G_{X\times X}^+]$. 

Thus, we have 

\proclaim{Proposition 2.20} The Hilbert series of the quadratic algebras 
$Q_X$,$A_X$ are equal to $(1-t)^{-N},(1-t)^{-N^2}$. 
\endproclaim

Thus, the Hilbert series are always as in the classical case $R=1$,
when we have usual polynomial algebras. 
 
According to \cite{FRT}, the algebra $A_X$ is a bialgebra, 
with coproduct and counit defined by
$$
\Delta(L)=L^{12}L^{13},\epsilon(L)=1.\tag 2.13
$$
Moreover, the algebra $Q_X$ is a left comodule over $A_X$, 
so that the coaction is an algebra homomorphism. 
This coaction is defined by 
$$
\Delta(q)=L^{12}q^{13},\tag 2.14
$$
and is the quantum analogue of the action 
of the monoid $\text{End}(V)$ on the space $V$.

Furthermore, the algebra $A_X$ can actually be extended 
to a Hopf algebra. Namely, define the algebra 
$\hat A_X$ to be generated by $A_X$ and 
an additional set of generators
$(L^{-1})_{xy}$, $x,y\in X$, 
with an additional defining relation 
$$
LL^{-1}=L^{-1}L=1,\tag 2.15
$$
where $L^{-1}:=\sum E_{xy}\o (L^{-1})_{xy}$. 

\proclaim{Proposition 2.21}
The algebra $\hat A_X$ is a Hopf algebra with the antipode defined by
$$
\gamma(L)=L^{-1},\gamma(L^{-1})=L.\tag 2.16
$$
\endproclaim

\demo{Proof} We need to check the antipode axiom and the fact 
that (2.16) extends to an antiautomorphism of $\hat A_X$. 
The first statement is immediate, and the second is checked by a direct 
computation. 
$\square$\enddemo

In particular, $\gamma^2=1$. 
Thus, $\hat A_X$ is an involutive Hopf algebra, which is a quantum 
analogue of the algebra of regular functions on the group $GL_N$. 

Now let $Comod(\hat A_X)$ be the category of finite-dimensional comodules
over $\hat A_X$. Since $\hat A_X$ is a Hopf algebra, this category 
is a rigid tensor category (for a definition of a rigid
tensor category, see e.g. \cite{DM}). If $R=1$, this category coincides with
the category of finite-dimensional representations of $GL_N$. 

\proclaim{Theorem 2.22} The category $Comod(\hat A_X)$ 
is equivalent to the category $Rep(GL_N)$ of finite-dimensional 
representations of $GL_N(\C)$ as a rigid tensor category.
\endproclaim

The proof of Theorem 2.22 is contained in the next section. 

\subhead 2.10. Proof of Theorem 2.22
\endsubhead

Let $\O(Mat_N)$ be the bialgebra of polynomial functions on $Mat_N(\C)$. 

\proclaim{Proposition 2.23} There exists 
a coalgebra isomorphism $\eta:\O(Mat_N)\to   
 A_X$. 
\endproclaim

\demo{Proof} Fix a labeling of elements of $X$ by 
indices $\{1,...,N\}$, so that $X=\{x_1,...,x_N\}$.  
We have $\O(Mat_N)=\C[T_{ij}]$, where $T=(T_{ij})$ is a matrix 
of indeterminates.
Define $\eta$ by the formula
$$
\eta(T_{i_1j_1}...T_{i_kj_k})=L_{y_1z_1}...L_{y_kz_k},\tag 2.17
$$
where $(y_1,...,y_k)=J_k^{-1}(x_{i_1},...,x_{i_k})$, 
$(z_1,...,z_k)=J_k^{-1}(x_{j_1},...,x_{j_k})$, where $J_k$ 
was defined in Section 1.3.  
By Proposition 1.7, the map $\eta$ is well defined, and is a linear 
isomorphism. It is easy to check that $\eta$ also respects the coproduct. 
Thus, $\eta$ is a coalgebra isomorphism. 
$\square$\enddemo

Proposition 2.23 implies that for any irreducible 
representation $W$ of $GL(V)$ ($V=\C^N$)
which occurs in $V^{\o k}$ for some $k$, we can define 
its direct image -- the corresponding comodule $\eta_*W$ of 
$A_X$. 

In particular, consider the 1-dimensional comodule
$Det=\eta_*(\Lambda^NV)$. It is a 1-dimensional space $\C v$, with 
coaction given by $\Delta(v)=D\o v$, where $D$ is an element of $A_X$. 
It is easy to see that this element is central and 
group-like (i.e. $\Delta(D)=D\o D$), and does not depend 
on the labeling of $X$. The element $D$ is called 
{\it the quantum determinant}.

It is clear that $D$ is not a zero divisor 
in $A_X$. 
Indeed, 
$D$ is invertible in $\hat A_X$, and 
$D^{-1}$ is the quantum determinant 
of $(L^{-1})^{t_1}$ (i.e. $D^{-1}$ is obtained from 
$(L^{-1})^{t_1}$ in the same way as $D$ is obtained from $L$). 
Therefore, it makes sense to consider 
the algebra $\tilde A_X=A_X\o_{\C[D]}\C[D,D^{-1}]$. 
This algebra inherits a bialgebra structure form $A_X$. 

We claim that the algebras $\hat A_X$, $\tilde A_X$ are isomorphic 
as bialgebras. Indeed, $L$ is invertible in $\tilde A_X$,
and  $L^{-1}=D^{-1}M_{N-1}(L)$, where $M_{N-1}(L)$ is 
a polynomial of $L$ of degree $N-1$ (the matrix of quantum minors). 
This allows to define an obvious isomorphism between the algebras
$\hat A_X$, $\tilde A_X$, which is clearly an isomorphism of bialgebras.

Let 
$\O(GL_N)$ be the Hopf algebra of polynomial functions on $GL_N(\C)$,

\proclaim{Proposition 2.24} There exists 
a coalgebra isomorphism $\hat\eta:\O(GL_N)\to   
 \hat A_X$. 
\endproclaim

\demo{Proof} It is easy to see that 
the map $\eta$ from Proposition 2.23 is $\C[D]$-linear
(for $\O(Mat_N)$, $D$ denotes the usual determinant).
Thus, $\hat\eta$ is obtained simply by tensoring $\eta$ over $\C[D]$
with $\C[D,D^{-1}]$.
$\square$\enddemo

Proposition 2.24 implies that the 
functor $\eta_*$ of direct image
induces an equivalence of abelian categories 
$Rep(GL_N)\to Comod(\hat A_R)$. So, to prove 
Theorem 2.22, it is enough to introduce a 
tensor structure on this functor. 

By definition, the tensor structure is 
a collection of functorial isomorphisms 
$J_{WU}:\eta_*(W)\o \eta_*(U)\to \eta_*(W\o U)$, satisfying 
the following compatibility condition:
$J_{W\o U,Y}(J_{WU}\o 1)=(1\o J_{UY})J_{W,U\o Y}$ \cite{DM}. 

Let us define $J_{WU}$ for $W=V^{\o m},U=V^{\o n}$. This should be an
operator $J_{V^{\o m}V^{\o n}}: V^{\o m+n}\to V^{\o m+n}$. We set
$$
J_{V^{\o m}V^{\o n}}=J_{m+n}(J_m^{-1}\o J_n^{-1}).\tag 2.18
$$

 It is easy to see that $J_{V^{\o m}V^{\o n}}$ is an intertwiner,
and that it satisfies the 2-cocycle condition. 
Further, by Proposition 1.7, $J_{V^{\o m}V^{\o n}}$
commutes with $S_m\times S_n$ (acting by permutations), 
so by Weyl duality it defines $J_{WU}$ for any irreducible representations
$W,U$ which occur in $V^{\o m}$ for some $m$. 
Now set $J_{W'U'}=J_{WU}$ if $W'=W\o (Det)^{\o k}$, $U'=U\o (Det)^l$. 
This defines $J_{WU}$ for any irreducible finite-dimensional 
representations $W,U$ of $GL_N$. Since the category of finite-dimensional
representations of $GL_N$ is semisimple, we have defined $J_{WU}$ for any
objects $W,U$ in $Rep(GL_N)$. We have

\proclaim{Proposition 2.25} The maps $J_{WU}$ define a tensor structure on 
$\eta_*$.
\endproclaim

\demo{Proof} Clear.
$\square$\enddemo

Theorem 2.22 is proved.

{\bf Remark.} The proof of Theorem 2.22 shows that the bialgebra 
$\hat A_X$ is obtained from $\O(GL_N)$ by twisting the multiplication,
in the sense of Drinfeld. 

\head 3. Methods of construction of nondegenerate symmetric sets.
\endhead

\subhead 3.1. Linear and affine solutions\endsubhead

In this section we  
will look for nondegenerate symmetric sets of the following form:
$X$ is an abelian group, and $S$ is an 
affine linear transformation of $X\times X$.
Such symmetric sets will be called {\it affine solutions}. 
Considering affine solutions was motivated by the results in \cite{Hi}. 

We will start with considering a special case, when
$S$ is an automorphism of $X\times X$. In this case, an affine solution 
will be called {\it a linear solution}. For a linear solution, 
$S$ has the form 
$$
S(x,y)=(ax+by,cx+dy),\ a,b,c,d\in \End X. \tag 3.1
$$
It is easy to check that for $S$ of the form (3.1) 
the braid relation is equivalent 
to the equations \cite{Hi}
$$
\gather
a(1-a)=bac,\ d(1-d)=cdb,\ ab=ba(1-d),\ ca=(1-d)ac,\ dc=cd(1-a),\\ 
bd=(1-a)db,\  cb-bc=ada-dad.\tag 3.2
\endgather
$$
It is also easy to see that the involutivity 
of $S$ is equivalent to the equations
$$
a^2+bc=1,\ cb+d^2=1,\ ab+bd=0,ca+dc=0.\tag 3.3
$$
Finally, the nondegeneracy condition is obviously equivalent 
to the condition that $b,c$ are invertible. 

\proclaim{Proposition 3.1} If $b,c$ are invertible, equations (3.2),(3.3) are 
equivalent to the equations 
$$
bab^{-1}=\frac{a}{a+1},c=b^{-1}(1-a^2),d=\frac{a}{a-1}.\tag 3.4
$$
\endproclaim

\demo{Proof} The second equation of (3.4) follows directly from (3.3). 
Also, (3.3) implies 
$$
a=-bdb^{-1}.\tag 3.5
$$
Therefore, multiplying the equation $bd=(1-a)db$ (which is in (3.2)) by $b^{-1}$ 
on the right, we get 
$$
-a=(1-a)d.\tag 3.6
$$
Since $b,c$ are invertible, so is $bc=1-a^2$, so 
$1-a$ is invertible. Thus, (3.6) implies the third equation of (3.4). 
Now the first equation of (3.4) follows from (3.5). 

Conversely, substituting (3.4) into (3.2),(3.3), it is easy to show 
by a direct calculation that they are identically satisfied. 
$\square$\enddemo

\proclaim{Corollary 3.2} 
A map $S$ of the form (3.1) is a linear solution 
if and only if $b,c$ are invertible, and (3.4) are satisfied. 
Thus, such solutions are in 1-1 correspondence with pairs $(a,b)$ such 
that $bab^{-1}=\frac{a}{a+1}$. 
\endproclaim

Now consider general affine solutions. Then $S$ has the form 
$$
S(x,y)=(ax+by+z,cx+dy+t), t,z\in X.\tag 3.7
$$
In this case, it is clear that the equations on $a,b,c,d$ are the same as 
before. The only equation for $z,t$ is obtained from 
the braid relation and has the form
$t=-b^{-1}(1+a)z$. Thus, we get

\proclaim{Proposition 3.3} 
A map $S$ of the form (3.7) is an affine solution 
if and only if $b,c$ are invertible, (3.4) are satisfied, and
$t=-b^{-1}(1+a)z$. 
Thus, such solutions are in 1-1 correspondence with triples (a,b,z) such 
that $bab^{-1}=\frac{a}{a+1}$. 
\endproclaim

Now consider examples of solutions of the equation 
$$
bab^{-1}=\frac{a}{a+1}.\tag 3.8
$$ 

{\bf Example 1.} \cite{Hi} Let $X=\Z/n \Z$. Then $\End X=\Z/n\Z$, 
which is commutative, so equation (3.8) reads $a=\frac{a}{a+1}$, which is 
equivalent to $a^2=0$ (and $b$ is any invertible element). 

{\bf Example 2.} Let $X=V^N$, where $V$ is an abelian group.
Then the algebra $Mat_N(\Z)$ of integer matrices is mapped into 
$\End X$. Thus, it is enough for us to construct a solution of 
$bab^{-1}=\frac{a}{a+1}$ in $Mat_N(\Z)$, such that $b\in GL_N(\Z)$. 
 
Let $a_{ij}=\delta_{i+1,j}$, and $b_{ij}=\left(\matrix j\\ i\endmatrix\right)$. 
Then $a,b$ satisfy (3.8). Indeed, this equation can be 
rewritten as $ab=ba+aba$, which at the level of matrix elements
reduces to the well-known identity
for binomial coefficients:
$$
\left(\matrix j\\ i+1\endmatrix\right) =
\left(\matrix j-1\\ i\endmatrix\right) +
\left(\matrix j-1\\ i+1\endmatrix\right). \tag 3.9
$$ 
We will use the following notation for this solution: 
$a=J_N$, $b=B_N$. 

In fact, all solutions of (3.8) in $Mat_N(\Z)$ can be obtained 
from $J_N,B_N$. Indeed, we have 

\proclaim{Lemma 3.4} Let $a,b$ be a solution of (3.8) in $Mat_N(\C)$. 
Then $a$ is nilpotent. 
\endproclaim

\demo{Proof} It follows from (3.8) that if $\l$ is an eigenvalue 
of $a$ then so is $\frac{\l}{\l+1}$. Therefore, if $\l\ne 0$, 
we get that $a$ has infinitely many distinct eigenvalues. 
This is impossible, so $\l=0$.
$\square$\enddemo

Thus, $a$ is nilpotent. Then,
by Jordan's theorem, $a$ can be reduced, over $\Bbb Q$, to Jordan normal form:
$a=J_{N_1}\oplus...\oplus J_{N_K}$, where 
$J_{N_l}\in Mat_{N_l}(\Z)$ are given by $(J_{N_l})_{ij}=\delta_{i+1,j}$. 
If $a$ is of this form, then $b=b_0A$, where $A$ commutes with $a$, and 
$b_0=B_{N_1}\oplus...\oplus B_{N_K}$.
Thus we have proved 

\proclaim{Proposition 3.5} Any solution of (3.8) in $Mat_N(\Z)$ with $b\in 
GL_N(\Z)$ is conjugate under $GL_N(\Bbb Q)$ to a solution of the form
$a=J_{N_1}\oplus...\oplus J_{N_K}$,
$b=(B_{N_1}\oplus...\oplus B_{N_K})A$, where $[A,a]=0$. 
\endproclaim

\proclaim{Proposition 3.6} If $V=\Z/p\Z$, where $p$ is a prime, 
and $N<p$, then any solution of (3.8) in $\End X$ is 
conjugate to a solution of the form 
given in Proposition 3.5. 
\endproclaim

\demo{Proof} Let $\l$ be an eigenvalue of $a$ over 
the algebraic closure $\overline{\Z/p\Z}$. 
Then $\frac{\l}{\l+1}$ is also an eigenvalue. Therefore, 
if $\l\ne 0$, $a$ has to have at least $p$ distinct eigenvalues.
Thus, $\l=0$ and hence $a$ is nilpotent.  
The rest of the proof is the same as for Proposition 3.5
(but instead of working over $\Bbb Q$ we work over $\Z/p\Z$). 
$\square$\enddemo

However, if $N\ge p$, other solutions are possible. 

{\bf Example: } $p=N=2$, $a=\left(\matrix 1&1\\ 1&0\endmatrix\right)$, 
$b=\left(\matrix 0&1\\ 1&0\endmatrix\right)$.
  
\subhead 3.2. Multipermutation solutions and equivariant 
fiber bundles\endsubhead

Let $(X,S)$ be a nondegenerate symmetric set. Then we can define 
another nondegenerate symmetric set $(\bar X,\bar S)$, such that 
there exists a surjective morphism $\mu:(X,S)\to (\bar X,\bar S)$.  
This is done as follows.

For $x,y\in X$ we will write $x\sim y$ if $f_x=f_y$. Since 
$g_x=T^{-1}f_x^{-1}T$, in this case we also have $g_x=g_y$. 
It is clear that $\sim$ is an equivalence relation. 
Let $\bar X=X/\sim$. For $x\in X$, let $\bar x$ denote the image of $x$ in
$\bar X$. 

By the definition, we have $x\sim y$ if and only if the images 
$\bar a_x,\bar a_y$ of $x,y$ in $G_X^0$ coincide. 
Thus, $\bar X$ is naturally identified with the image 
of the map $\bar a:X\to G_X^0$ given by $\bar a_x=f_x^{-1}$. 
  
Let $A$ be the abelian group defined in Section 2.5.
It is easy to see that $x\sim y$ if and only if 
$\bar t_x=\bar t_y$, where $\bar t_x$ is the image of $t_x$ 
under the homomorphism $\Z^X\to A$. 
This follows from the fact that $\bar a_x=\bar\pi^{-1}(\bar t_x)$, where  
$\bar\pi:G_X^0\to A$ is the bijective 1-cocycle defined 
in Section 2.5. Thus, $\bar X$ is naturally identified 
with the image of $X\subset \Z^X$ under the map $\Z^X\to A$. 

The subset $\bar X\subset A$ is invariant under the action of 
$G_X$, i.e. under operators $f_x^{-1}$. Therefore, it is also 
invariant under $g_x$. Thus,
according to the remark in Section 2.4, 
$\bar X$ has a structure of a nondegenerate 
 symmetric set, with 
$\bar S(\bar x,\bar y)=(\overline{g_x(y)},\overline{f_y(x)})$, where
$\bar x,\bar y\in \bar X$, and $x,y$ are any preimages of 
$\bar x,\bar y$ in $X$. We will call $(\bar X,\bar S)$ the 
{\it retraction} of $(X,S)$, and denote it by $Ret(X,S)$. 
 
\proclaim{Definition 3.1} A solution 
$(X,S)$ will be called a multipermutation solution of level $n$ if 
$n$ is the minimal nonnegative integer such that
$|Ret^n(X,S)|=1$ (i.e. $Ret^n(X,S)$ is finite of size 1). 

A solution $(X,S)$ is called irretractable if $Ret(X,S)=(X,S)$
(i.e. if $\sim$ is a trivial equivalence relation).
\endproclaim

In particular, a multipermutation solution of level 0 is the trivial
solution for $|X|=1$, and a multipermutation solution of level 1 is a
permutation solution. On the other hand, as we will see below, there
exist irretractable affine solutions for $|X|=4$.

In terms of the bijective cocycle datum $(G_X^0,A,\rho,\pi)$
corresponding to $(X,S)$ in the manner of Section 2.5, there is a nice
formulation of retractions.  Indeed, we have seen that $\bar X$ is
just the image of $X$ in $A$, so the retraction is just given by a
set-structure on $(G_X^0,A,\rho,\pi)$ given by $(\bar X,\rho,\phi)$
where $\phi$ is the tautological embedding $\bar X \rightarrow A$, and
the action $\rho$ of $G$ on $A$ defines its action on $X$.
We see that the unique datum with a
faithful generating set-structure for $(X,S)$ is given by $\bar
G=G_X^0/H$, $\bar A=A/K$ where $H=Ker(\rho)$, $K=\pi(H)$ along with
the maps $\bar \rho,
\bar \pi$ descending from $\rho, \pi)$.  Then, the faithful generating set-structure
for $(\bar G, \bar A, \bar \rho, \bar \pi)$ is given by $(X, \rho,
\bar \phi)$ where $\rho$ is the action of $G$ on $X$ as before, 
which descends modulo $H$ to an action of $\bar G$ on $X$, and the map
$\bar \phi: \bar X \rightarrow
\bar A$ follows from the natural map $A \rightarrow A/K$. 

Thus, given a datum $(G,A,\rho,\pi)$ define its retraction by
$Ret(G,A,\rho,\pi) =$ \linebreak $(G/H,A/K,\bar \rho, \bar \pi)$ as
above.  We have seen that this definition corresponds to that for
nondegenerate symmetric sets in the case of faithful generating
set-structures.

In view of this we have the following proposition:

\proclaim{Proposition}  
If the group $A$ obtained from a nondegenerate symmetric set $(X,S)$ as in Section 2.5
is finite and cyclic, then $(X,S)$ is a multipermutation solution.
\endproclaim

\demo{Proof} 
Assume $|A|>1$. Then, since $|Aut(A)|=\phi(|A|) < |A|$, where $\phi$
is the Euler $\phi$-function, $|Ker(\rho)|>1$ and so
$(G_X^0,A,\rho,\pi)$ is retractable yielding a new datum with a
smaller cyclic group $\bar A$.  Inductively we find that for some $n$,
$Ret^n(G_X^0,A,\rho,\pi)$ $=(G',A',\rho',\pi')$ is trivial
($A'=G'=\{0\}$) since $|A|$ was finite; clearly $n=0$ if $|A|=1$. This
means that $Ret^n(X,S)$ is given by the trivial datum together with
any faithful generating set-structure, which we find is just the
trivial solution on $|X|$ elements.  Hence $|Ret^{n+1}(X,S)|=1$.
$\square$\enddemo

%In fact, any other group $A'$ (other than $A/K$) that
%satisfies this property with $\rho'$ the required action on $X'$ is
%isomorphic to $A/K$.  This is because, given any group $A'$ satisfying
%this property, the action $\rho'$ of $G'$ on $X' \subset G'$
%determines the action $\rho'$ of $G'$ on all of $G'$ inductively as
%follows: Let $x, y \in G'$, $z \in X$, such that we already know the
%action of $G'$ on $y$.  Then $x * yz = \pi^{-1}(x * (\pi(z) +
%\pi(z^{-1} * y))) = \pi^{-1}(\pi(x * z) + \pi(xz^{-1} * y)) = ((x *
%z)xz^{-1} * y)(x * z)$ and similarly $x * yz^{-1} = \pi^{-1}(x *
%(\pi(z^{-1}) + \pi(z * y))) = \pi^{-1}(-(xz * z) + xz * y) =
%(T^{-1}(xz * z)^{-1}xz * y)(T^{-1}(xz * z)^{-1})$ with $T(z)=z * z$ as
%usual is a bijection $X' \rightarrow X'$.  Then this uniquely determines addition in $A$ by $\pi(x)+\pi(y)=\pi((x * y)y)$.

%this determines the action of $G'$ on all of $G'$ inductively, given that $G'$ is generated by $X'$:  $x * yz = x * (

Now we will consider solutions $(Y,S_Y)$ such that 
$Ret(Y,S_Y)$ is a fixed solution $(X,S)$. 
Such a solution $(Y,S_Y)$ will be called a blow-up of $(X,S)$. 
   
Let $BL(X,S)$ denote the category of all 
such solutions, where morphisms are homomorphisms of solutions
which become identity under retraction.

Our goal is to describe this category in group-theoretical terms.
Recall the following standard definition.

\proclaim{Definition 3.2} Let $X$ be a set and $G$ a group acting on $X$. 
A $G$-equivariant fiber bundle over $X$ is a set $Y$ equipped with a 
surjective map $p:Y\to X$, and an action $\rho$ 
of $G$ on $Y$ which respects $p$ 
and descends under $p$ to the action of $G$ on $X$.
\endproclaim

Denote by $Bun(X,G)$ the category of 
 $G$-equivariant fiber bundles over $X$,
where morphisms are $G$-invariant bundle mappings. 
Denote by $Bun_f(X,G_X)$ the full subcategory of $Bun(X,G_X)$ 
which consists of such bundles that $\rho(x)\ne \rho(y)$ for 
$x,y\in X\subset G_X$, $x\ne y$. Objects of $Bun_f$ will 
be called {\it faithful} bundles.  

\proclaim{Theorem 3.7} The categories $BL(X,S)$ and 
$Bun_f(X,G_X)$ are equivalent. In particular, there is a 1-1 correspondence 
between isomorphism classes of blow-ups of size $n$, and isomorphism  
classes of faithful bundles of size $n$.
\endproclaim

\demo{Proof} To prove the theorem, it is enough to construct 
two functors, $E: BL(X,S)\to Bun_f(X,G_X)$, and $E':Bun_f(X,G_X)\to 
BL(X,S)$, so that $E\circ E'=id,E'\circ E=id$. 

Let us construct $E$. Let $(Y,S_Y)\in BL(X,S)$. Then by definition
$X=\bar Y$ and thus we have a natural surjective map 
$p:Y\to X$. Moreover, the group $G_X$ acts on $Y$ by 
$\rho(x)(y)=f_{\tilde x}^{-1}(y)$, where $\tilde x\in Y$ is any 
lifting of $x$ to $Y$. As we explained, the action $\rho$ 
respects the map $p$ and  descends to the standard action of $G_X$ on $X$ 
under $p$. Thus, $(Y,p,\rho)\in Bun(X,G_X)$. 
It is easy to see that in fact $(Y,p,\rho)\in Bun_f(X,G_X)$.
Set $E(Y,S_Y)=(Y,p,\rho)$. 
It is clear that $E$ is a functor. 

Now let us construct $E'$. Let $(Y,p,\rho)\in Bun_f(X,G_X)$. 
For $y\in Y$, define $f_y:Y\to Y$ by $f_y^{-1}=\rho(p(y))$. 
Define $g_y:Y\to Y$ by the usual formula $g_y(z)=f^{-1}_{f_z(y)}(z)$. 
Set $S_Y(y,z)=(g_y(z),f_z(y))$. Then $(Y,S_Y)$ is involutive. 

Let us show that $(Y,S_Y)$ is nondegenerate. It is clear that $f_y$ is 
invertible. To show that $g_y$ is invertible, it is enough to show that 
$T$ is invertible, due to Proposition 2.2(a). Recall that $T(y)=f_y^{-1}(y)$. 
To show that $T$ is invertible, it is enough to show that the equation 
$f_y^{-1}(y)=z$ has a unique solution for any $z$. 

Let $z\in Y$ and $\bar z\in X$ be its equivalence class. Since the map 
$T$ for $X$ is invertible, we can find a unique $\bar y\in X$ such that 
$f_{\bar y}^{-1}(\bar y)=\bar z$. Now to solve the equation $T(y)=z$, 
 is the same as to find an element $y\in p^{-1}(\bar y)$ such that 
$f_{\bar y}^{-1}(y)=z$. Such an element exists and unique, since by
definition 
$f_{\bar y}^{-1}$ induces an isomorphism of the fibers  
 $p^{-1}(\bar z)$ and $p^{-1}(\bar y)$.  
Thus, $(Y,S_Y)$ is nondegenerate. 

Since $\rho$ is an action of $G_X$ on $Y$, the maps $f_y$ satisfy 
condition (i) of Proposition 2.1. Therefore, by Proposition 2.2, 
$(Y,S_Y)$ is a nondegenerate symmetric set, and 
$Ret(Y,S_Y)=(X,S)$, since $Y\in Bun_f$. Set $E'(Y,p,\rho)=(Y,S_Y)$. 
Clearly, $E'$ is a functor. 

The identities $E\circ E'=id, E'\circ E=id$ are obvious. 
The theorem is proved.
$\square$\enddemo

Let $K_x$ denote the stabilizer of a point $x\in X$ in $G_X$.
 
\proclaim{Proposition 3.8} Let $(Y,S_Y)\in BL(X,S)$. Then $Y$ is 
indecomposable if and only if $X$ is indecomposable, and 
$K_x$ acts transitively on $p^{-1}(x)$ for some $x\in X$. 
\endproclaim

\demo{Proof} Clear.
$\square$\enddemo

{\bf Remark.} Theorem 3.7 and Proposition 3.8 reduce the classification 
of indecomposable solutions to classification  
of irretractable indecomposable solutions, 
modulo the group-theoretical question of classification 
of faithful bundles. 

 From now on and till the end of Chapter 3 we assume that the set $X$ is
finite.
\subhead 3.3. Indecomposable multipermutation solutions of level 2
\endsubhead

Let $(Y,S_Y)$ be an indecomposable multipermutation solution of level 2.
Then $Ret(Y,S_Y)=(X,S)$, where $X=\Z/m\Z$, $m>1$, and $S(x,y)=(y-1,x+1)$.

Let us write $Y$ as a direct product $Y=X\times Z$, where $Z=\{1,...,n\}$. 
We may suppose that the natural map  $p:Y\to X$ is the projection to the 
first component. 

Now consider the group $G_X$. This group is the subgroup 
of $Aut(X)\ltimes \Z^X$ consisting of the elements
$c^k(b_0,...,b_{m-1})$, where $c:X\to X$ is given by 
$c(x)=x-1$, and $b_i\in \Z$ are such that $\sum b_i=k\text{ mod }m$. 
In particular, the group $K_x$ of elements stabilizing the point $x$ 
is the group $\Gamma$ of vectors  
$(b_0,...,b_{m-1})\in \Z^X$ such that $\sum b_i$ is divisible by $m$.
(This is the same group $\Gamma$ as we considered before).

As we know, blow-ups of $X$ correspond to faithful 
$G_X$-equivariant fiber bundles on $X$. It is clear that 
equivariant fiber bundles on $X$ correspond to 
$\Gamma$-actions on $Z$, via $Z\to G_X\times_{\Gamma}Z$. 

Thus, it remains to classify transitive actions of $\Gamma$ on 
finite sets $Z$ of size $n$. Since $\Gamma$ is a free abelian group, such 
actions correspond to sublattices $L$ in $\Gamma$ of index $n$
($Z=\Gamma/L$). When all such sublattices are found, one should  
separate those which define faithful bundles. 

{\bf Example 1.} $m=2,n=2$. In this case we have to classify 
transitive actions of $\Gamma$ on a 2-element set, i.e. surjective 
maps $\Gamma\to \Z/2\Z$. Since $\Gamma$ is the sublattice in 
$\Z^2$ generated by $(1,1),(1,-1)$, we have 3 choices:

1. $(1,1)\to 1,(1,-1)\to 1$.
 
2. $(1,1)\to 0,(1,-1)\to 1$.

3. $(1,1)\to 1,(1,-1)\to 0$.

One can check that the first two choices define faithful bundles, 
while the third choice does not. 
Thus, there are two 
indecomposable multipermutation solutions of level 2 for $|Y|=4$
(one can show that choices 1 and 2 define non-equivalent solutions).
The total number of indecomposable 
multipermutation solutions for $|Y|=4$ is three, since 
we also have the cyclic permutation solution. 
It turns out, however, that there are two more indecomposable 
solutions for $|X|=4$ -- 
they are irretractable affine solutions (see below).

{\bf Example 2.} $m=3,n=2$. In this case we have to classify 
transitive actions of $\Gamma$ on a 2-element set, i.e. surjective 
maps $\Gamma\to \Z/2\Z$. Since $\Gamma$ is the sublattice in 
$\Z^3$ spanned by $(1,0,-1),(1,-1,0),(1,1,1)$, we have 7 choices:

1. $(1,0,-1)\to 1,(1,-1,0)\to 0,(1,1,1)\to 0$.

2. $(1,0,-1)\to 0,(1,-1,0)\to 1,(1,1,1)\to 0$.

3. $(1,0,-1)\to 0,(1,-1,0)\to 0,(1,1,1)\to 1$.

4. $(1,0,-1)\to 1,(1,-1,0)\to 1,(1,1,1)\to 0$.

5. $(1,0,-1)\to 1,(1,-1,0)\to 0,(1,1,1)\to 1$.

6. $(1,0,-1)\to 0,(1,-1,0)\to 1,(1,1,1)\to 1$.

7. $(1,0,-1)\to 1,(1,-1,0)\to 1,(1,1,1)\to 1$.

All of these choices but choice 3 define faithful bundles, which 
give rise to non-equivalent solutions. Thus, we get 6 indecomposable 
solutions for $|Y|=6$. 

{\bf Example 3.} $m=2,n=3$. In this case we have to classify 
transitive actions of $\Gamma$ on a 3-element set, i.e. surjective 
maps $\Gamma\to \Z/3\Z$. Since $\Gamma$ is the sublattice in 
$\Z^2$ generated by $(1,1),(1,-1)$, we have 8 choices:

1. $(1,1)\to 1,(1,-1)\to 1$.
 
2. $(1,1)\to 0,(1,-1)\to 1$.

3. $(1,1)\to 1,(1,-1)\to 0$.

4. $(1,1)\to 2,(1,-1)\to 2$.
 
5. $(1,1)\to 0,(1,-1)\to 2$.

6. $(1,1)\to 2,(1,-1)\to 0$.

7. $(1,1)\to 2,(1,-1)\to 1$.

8. $(1,1)\to 1,(1,-1)\to 2$.

All choices except 3 and 6 define faithful bundles, 
and we have the following isomorphisms between corresponding 
solutions: 1-4,2-5,7-8. Thus, we get 3 more indecomposable solutions
for $|Y|=6$. 

The total number of indecomposable multipermutation 
solutions 
for $|X|=6$ is 6+3+1=10, since we also have the cyclic permutation solution. 
A computer calculation shows (see below) that these are all indecomposable 
solutions for $|Y|=6$. 

\subhead 3.4. Twisted unions and generalized twisted unions\endsubhead

Let $X,Y$ be finite nondegenerate symmetric sets, and 
$Z=X\cup Y$ be their union. 
 
\proclaim{Definition 3.3} $Z$ is called a twisted 
union of $X$ and $Y$ if the map $S_Z:X\times Y\to Y\times X$ 
is given by the formula 
$$
S_Z(x,y)=(g(y),f(x)),\tag 3.10
$$ where $g:Y\to Y$, $f:X\to X$ are permutations. 
\endproclaim

It follows from involutivity that for a twisted union, 
$$
S_Z(y,x)=(f^{-1}(x),g^{-1}(y)).\tag 3.11
$$

The classification of twisted unions is very simple. 

\proclaim{Proposition 3.9} Formulas (3.10),(3.11) define a union of $X$ and $Y$
if and only if $f$ preserves $S_X$ and $g$ preserves $S_Y$. 
\endproclaim

\demo{Proof} The nondegeneracy and involutivity of $S_Z$ are automatic, 
and the braid relation is easily shown to be equivalent to 
the condition that $f$ preserves $S_X$ and $g$ preserves $S_Y$. 
$\square$\enddemo

{\bf Example 1.} Any permutation solution is naturally a twisted union 
of cyclic permutation solutions. 

{\bf Example 2.} Let $|X|=1$, and $Z=X\cup Y$ be a union. 
Then $Z$ is obviously a twisted union.   

At this point, it is easy for us to classify all solutions with 
$|X|\le 3$. Indeed, any such indecomposable solution is a cyclic
permutation solution by Theorem 2.12, and any decomposable one is a twisted 
union. However, already for $|X|=4$ there are unions which are not twisted 
unions. This encourages one to introduce the notion of a generalized 
twisted union. 

\proclaim{Definition 3.4} A union $Z$ of solutions $X$ and $Y$ 
is called a generalized twisted 
union of $X$ and $Y$ if the map $S_Z:X\times Y\to Y\times X$ 
is given by the formula 
$$
S_Z(x,y)=(g_x(y),f_y(x)),\tag 3.12
$$ 
where $g_x:Y\to Y$, $f_y:X\to X$ are permutations, such that 
the permutation $g_{f_y(x)}:Y\to Y$, $x\in X$, is independent of $y\in Y$, and 
the permutation $f_{g_x(y)}:X\to X$, $y\in Y$, is independent of $x\in X$. 
\endproclaim

For a generalized twisted union, we will write 
$g_{f_y(x)}$ and $f_{g_x(y)}$ as 
$g_{f_*(x)}$ and $f_{g_*(y)}$ (for $x\in X$, $y\in Y$).

It is easy to check that in a generalized twisted union, 
the permutation $g_{f_y^{-1}(x)}$ does not depend on $y$, and 
$f_{g_x^{-1}(y)}$ does not depend on $x$. Thus, we will denote 
these permutations by $g_{f_*^{-1}(x)}$ and $f_{g_*^{-1}(y)}$.
 
It follows from involutivity that in a generalized twisted 
union, 
$$
S_Z(y,x)=(f^{-1}_{g_*^{-1}(y)}(x),g^{-1}_{f_*^{-1}(x)}(y)),\tag 3.13
$$

It is clear that a twisted union is a special case of a generalized 
twisted union. 

\proclaim{Proposition 3.10} Formulas (3.12),(3.13) define a generalized 
twisted union if and only if the following conditions 
are simultaneously satisfied: 

(i) The assignments $y\to f_y^{-1}$, $x\to g_x$ define left actions 
of $G_Y$ on $X$ and of $G_X$ on $Y$. 

(ii) The map $S_X$ commutes with $f_{g_*(y)}\times f_y$. 
The map $S_Y$ commutes with $g_x\times g_{f_*(x)}$. 
\endproclaim

\demo{Proof} As in Proposition 3.9, involutivity and nondegeneracy are 
automatic. It is easy to check that the braid relation 
is equivalent to conditions (i) and (ii). 
$\square$\enddemo

{\bf Remark.} Note that for a generalized twisted union, 
the assignments $y\to f^{-1}_{g_*^{-1}(y)}$ and 
$x\to g_{f_*^{-1}(x)}$ are actions of $G_Y$ on $X$ and $G_X$ on $Y$. 
Let us call them the modified actions, as opposed to the standard 
actions $y\to f_y^{-1}$, $x\to g_x$. 
Thus, condition (ii) reads that $S_X$
is invariant under the product of the standard and the modified action, 
and similarly for $S_Y$. 

{\bf Example.} Any multipermutation solution of level 2 is a generalized 
twisted union of indecomposable multipermutation solutions of level $\le 2$.

\subhead 3.5. Solutions for $|X| \le 8$\endsubhead

A computer program in C generated all solutions $(S,X)$ up to
isomorphism.  Then, programs in Perl classified the solutions. Below
we summarize the results in a table.

The following abbreviations will be used: s. = solutions; t.u. =
twisted unions; g.t.u. = generalized twisted unions; id. =
indecomposable; d. = decomposable; ir. = irretractable; a. = affine;
mp. = multipermutation.  The table gives the number of distinct maps
up to isomorphism for $|X| \leq 8$:
\vskip 12 pt

$$
\boxed{\matrix 
|X| & s. & d.s. & t.u. & g.t.u. &
id. s. & id. mp. s. & id. ir. s. & id. ir. a. \\ 
1 & 1 & 0 & 0 & 0 & 1 & 1 & 0 & 0 \\ 
2 & 2 & 1 & 1 & 1 & 1 & 1 & 0 & 0 \\ 
3 & 5 & 4 & 4 & 4 & 1 & 1 & 0 & 0\\ 
4 & 23 & 18 & 16 & 18 & 5 & 3 & 2 & 2 \\ 
5 & 88 & 87 & 84 & 87 & 1 & 1 & 0 & 0 \\ 
6 & 595 & 585 & 425 & 585 & 10 & 10 & 0 & 0 \\
7 & 3456 & 3455 & 3270 & 3455 & 1 & 1 & 0 & 0 \\
8 & 34528 & 34430 & 23856 & 34350 & 98 & 37 & 47 & 0  \endmatrix }
$$

\vskip 12 pt 

For $|X| \leq 7$, all decomposable solutions were found to be
generalized twisted unions, and all indecomposables except for two
affine solutions turned out to be multipermutation solutions.

However, for $|X|=8$, 47 solutions were found to be irretractable
indecomposables, none of which are affine.  The 14 solutions that are
not irretractable and not multipermutation are clearly blow-ups of
the two irretractable indecomposable solutions for $|X|=4$.
Furthermore, 80 decomposable solutions for $|X|=8$ are not generalized
twisted unions.

\head 4. Power series solutions.\endhead

\subhead 4.1. The definition of a power series solution
\endsubhead
 
The notion of a linear solution, introduced in Section 3.1, 
can be generalized, by defining the notion 
of a power series solution. This is done as follows. 

Let $K$ be a ring, and $D^N$ be the formal N-dimensional polydisk 
over $K$. A formal mapping $D^N\to D^M$ is, by definition, 
a vector $\phi=(\phi_1,...,\phi_M)$, where 
$\phi_i\in K[[x_1,...,x_N]]$ are power series with zero free term. 

We should remember that $D^N$ is not a set but a formal scheme, and 
thus $\phi$ is not a mapping in the usual sense. 
To pass to sets and mappings, let $I$ be 
a nilpotent commutative algebra over $K$.
This means, $I$ is an algebra over $K$ (without unit), and 
any element of $I$ is nilpotent. In this case, 
any formal power series $\psi\in K[[x_1,...,x_N]]$ with zero free term 
defines a mapping $I^N\to I$.
Thus, a formal mapping $D^N\to D^M$ defines a mapping $I^N\to I^M$. 

Let $X=I^N$. Let $S:X^2\to X^2$ be a mapping, such that 
$(X,S)$ is a nondegenerate symmetric set. 

\proclaim{Definition 4.1} We will say that $(X,S)$ is a power series solution
if $S$ is induced by a formal mapping $(D^N)^2\to (D^N)^2$. 
\endproclaim

It is clear that a linear solution is a special case of a power 
series solution. Indeed, let $I$ be any abelian group. 
Equip $I$ with a $\Z$-algebra structure by defining its multiplication to be 
zero. Power series solutions for such $I$ are the same thing as linear 
solutions, since all nonlinear terms in the power series are automatically 
zero. 

Let $x,y$ be N-dimensional 
vectors of indeterminates, and $S:(D^N)^2\to (D^N)^2$ be a formal mapping.
The mapping $S$ is given by $S(x,y)=(S_1(x,y),S_2(x,y))$, where
$S_i(x,y)$ are formal mappings  $(D^N)^2\to D^N$. 

\proclaim{Definition 4.2} $S$ is called a universal power series 
solution over $K$ if for any nilpotent $K$-algebra $I$
the series $S$ defines the structure of a nondegenerate 
symmetric set on $X=I^N$. 
\endproclaim

We have the following simple proposition. 

\proclaim{Proposition 4.1} $S$ is a formal power series solution if and only if
it satisfies equations (1.1),(1.2), and 
the nondegeneracy condition: the series $S_1(x,*)$, $S_2(*,y)$ and invertible 
for fixed $x,y$. 
\endproclaim

\demo{Proof} The ``if'' statement is obvious. 
The ``only if'' statement: set $I=J/J^m$, where 
$J$ is the ideal in $K[[z_1,...,z_l]]$
consisting of series with the zero free term, and
$m,l$ are arbitrary integers. 
Then the condition that $(X,S)$ is a solution 
easily implies the claim. 
$\square$\enddemo

\proclaim{Definition 4.2} Two universal power series solutions 
$S,S'$ are said to be isomorphic if there exists 
an invertible formal mapping $\phi:D^N\to D^N$ 
such that $S'=(\phi\times\phi)S(\phi^{-1}\times \phi^{-1})$. 
\endproclaim

 From now on we will be interested only in universal power series solutions,
 and drop the word ``universal'' in our discussions. 
For simplicity we will assume that $K$ is a field.    
 
\subhead 4.2. Permutation power series solutions
\endsubhead

 We are interested in classification of power series solutions, up to 
isomorphism. 

Consider first the special case of this problem: classification of 
permutation solutions. By definition, a permutation solution 
has the form $S(x,y)=(f(y),f^{-1}(x))$, where $f:D^N\to D^N$ is an 
invertible formal mapping. Two such solutions corresponding to 
formal mappings $f,f':D^N\to D^N$ are isomorphic if there 
exists an invertible formal mapping $\phi$ such that $f'=\phi f\phi^{-1}$. 
Thus, the problem of classifying permutation power series solutions is 
equivalent to classifying conjugacy classes in the group $Diff(D^N)$
of formal diffeomorphisms (=invertible formal mappings) of $D^N$ into itself. 

It is well known that this problem is ``wild'' (i.e. impossible 
to solve effectively) for $N>1$. On the other hand, for $N=1$, 
it is ``tame'', and the solution is well known. It is given
by the following proposition (see \cite{Ar}). 

\proclaim{Proposition 4.2} Let $f(x)=\sum_{m\ge 1}a_mx^m$, $a_1\ne 0$. 

(i) If $a_1$ is not a root of unity, then $f(x)$ is conjugate 
to the linear map $f(x)=a_1x$. 

(ii)If $a_1$ is a root of unity of order $k$, 
then $f(x)$ is conjugate to normal form (i), or
to a map $f(x)=a_1x+x^{kr+1}+cx^{2kr+1}$, $r\ge 1$,  
and the parameters
$r,c$ are completely determined by $f$.

(iii) Normal forms (i),(ii) are never conjugate to each other. 
\endproclaim

Part (i) of Proposition 4.2 has 
the following multivariable generalization, due to 
Poincare (cf \cite{Ar}). 

We will say that a set of complex numbers $\{\l_1,...,\l_N\}$
is resonance free if the equation $\l_i=\l_{i_1}...\l_{i_m}$ 
is not satisfied for $m\ge 2$ and $1\le i,i_1,...,i_m\le N$. 
For example, the set $\{\l\}$ is resonance free if and only if 
$\l$ is not a root of unity.   

\proclaim{Proposition 4.3} Let $f(x)$ be a formal diffeomorphism of $D^N$, 
whose linear part has resonance free eigenvalues. Then $f(x)$ is conjugate 
to its linear part. 
\endproclaim
 
\subhead 4.3. Linearization of power series solutions 
\endsubhead

Consider a power series solution 
given by a formal
diffeomorphism \linebreak $R:D^N\times D^N\to D^N\times D^N$.
It is clear that the linear part $R_1$ of $R$ (i.e. 
the collection of all linear terms of $R$) defines a linear solution.
Proposition 4.3 gives rise to the question: when is a power series solution
isomorphic to its linear part? In this section we will partially 
answer this question. 

We will first consider the case when the linear part of $R(x,y)$ 
is the permutation solution $R_1(x,y)=(b^{-1}x,by)$. 
Suppose that the lowest degree of terms in the series $R-R_1$ is
$m$. 
Then we can write 
$$
R(x,y)=(b^{-1}x+P(x,y),by+Q(x,y))\text{ mod degree m+1},\tag 4.1
$$
where $P,Q$ are polynomials of degree $m$. 
The unitarity condition for $R$ gives the 
following equation for $P,Q$:
$$
P(x,y)=-b^{-1}Q(by,b^{-1}x).\tag 4.2
$$
Given (4.2), 
the quantum Yang-Baxter equation gives the following equation for $Q$:
$$
bQ(b^{-1}x,z)-Q(x,bz)=bQ(y,z)-Q(by,bz).\tag 4.3
$$

Introduce a new function $U(y,z)=bQ(y,z)-Q(by,bz)$. 
Then equation (4.3) can be rewritten in the form
$U(b^{-1}x,z)=U(y,z)$, which implies that $U(y,z)$ is independent of
$y$: $U(y,z)=U(z)$. Thus, we obtain the equation
$$
bQ(y,z)-Q(by,bz)=U(z).\tag 4.4
$$

\proclaim{Proposition 4.4} If $b$ has resonance free eigenvalues, 
then $Q(y,z)$ is independent of $y$.
\endproclaim

\demo{Proof} 

Let $g\in GL_N(K)$, $h\in GL_M(K)$, and 
$H$ the space of formal mappings $D^M\to D^N$. 
Let $A:H\to H$ be the linear operator defined by 
$(Af)(x)=gf(x)-f(hx)$. Let 
$\l_i$ be the eigenvalues of $g$, $\mu_j$ the eigenvalues of $h$.

{\bf Lemma.} If $\l_i\ne \mu_{j_1}...\mu_{j_m}$ for any $m\ge 2$
then the operator $A$ is invertible.

{\it Proof of the Lemma.} Let $H_m\subset H$ be the space of polynomials 
of degree $m$. It is clear that the eigenvalues of $A$ on $H_m$ 
are $\l_i-\mu_{i_1}...\mu_{i_m}$, so they are not zero.
The Lemma is proved.

Now take $M=2N$, $g=b$, $h=b\oplus b$. Then the Lemma implies 
that if the eigenvalues $\l_i$ of $b$ are resonance free,
the operator $Q(y,z)\to bQ(y,z)-Q(by,bz)$ is invertible. 
Therefore, equation (4.4) for a fixed $U(z)$ has a unique solution. 
This solution is obviously $y$-independent, since it does not change under 
rescalings of $y$. The Proposition is proved. 
$\square$\enddemo

\proclaim{Corollary 4.5} If eigenvalues of $b$ are resonance free
then any power series solution solution $R(x,y)$ whose linear part  
is $(bx,b^{-1}y)$ is isomorphic to its linear part
(i.e. is a permutation solution with the same linear part). 
\endproclaim

\demo{Proof} We will prove the statement modulo terms 
of degree $m+1$ by induction in $m$
Modulo terms of degree 2 the statement is a tautology.
Suppose we know it modulo degree $m$, and want to prove 
it modulo degree m+1. By the induction assumption, $R$ has 
the form (4.1). We proved that in this case $Q(x,y)$ is independent
on $x$ and thus by (4.2) $P(x,y)$ is independent on $y$.
Thus, modulo terms of degree $m+1$ the solution $R$ is a permutation 
solution, and thus by Poincare theorem (Proposition 4.3) it can be linearized
modulo degree $m+1$. The Corollary is proved.
$\square$\enddemo

Now we will generalize Corollary 4.5 to 
the case of an arbitrary linear part
$R_1(x,y)=(cx+dy,ax+by)$. For simplicity 
we will assume that $char K=0$ or that $N<char K$. 

\proclaim{Proposition 4.6} If the eigenvalues of $b$ are resonance free, 
then any solution $R$ with linear part 
$R_1(x,y)=(cx+dy,ax+by)$ is isomorphic to this linear part. 
\endproclaim

\demo{Proof}
In our case, the classification of linear solutions is as in 
Proposition 3.5. In particular, 
we have the following important property: there exists a basis 
of the N-dimensional space
(over the algebraic closure $\bar K$) in which the matrices
$b,c$ are upper triangular, and $a,d$ strictly upper triangular. 

Now consider $R(x,y)$ of the form
$$
R(x,y)=(cx+dy+P(x,y),ax+by+Q(x,y))\text{ mod degree m+1},\tag 4.5
$$
where $P,Q$ are polynomials of degree $m$. 
The unitarity condition for $R$ gives the 
following equation for $P,Q$:
$$
P(x,y)=-b^{-1}aQ(x,y)-b^{-1}Q(by,b^{-1}x).\tag 4.6
$$
Given (4.6), 
the quantum Yang-Baxter equation gives the following equation for $Q$:
$$
\gather
a(a+1)^{-1}Q(x,y)-a(a+1)^{-1}Q(ax+by,cx+dy)+bQ(cx+dy)=\\
bQ(y,z)+Q(x,ay+bz)-Q(ax+by,acx+ady+bz).\tag 4.7
\endgather
$$

{\bf Lemma.} If $b$ has resonance free eigenvalues, then any solution
of equation (4.7) has the form  
$$
Q(x,y)=aT(x)+bT(y)-T(ax+by),\tag 4.8
$$
where $T$ is a suitable polynomial of degree $m$. 

{\it Proof of the Lemma.}  First of all, (4.8) is actually a solution of 
(4.7). Indeed, if we conjugate the linear solution $R_1$ by the permutation 
$x\to x-T(x)$, we will get exactly (4.8). 
 
Denote by $D(a,b)$ the dimension of the space of solutions 
of equation (4.7). Since the eigenvalues of $b$ are resonance free, 
the assignment $T(x)\to aT(x)+bT(y)-T(ax+by)$ is injective. 
Thus, $D(a,b)\ge N\left(\matrix N+m-1\\ m\endmatrix\right)$
(the dimension of the space of polynomials $T$). 

We will work in the aforementioned basis, in which 
$a,b,c,d$ are upper triangular. 
Let $b_{diag}$ be the diagonal part of $b$ in this basis. 
We showed above that 
the lemma is true for $a=0$, i.e.
$D(0,b_{diag})=N\left(\matrix N+m-1\\ m\endmatrix\right)$.
So it is enough to show that $D(a,b)\le D(0,b_{diag})$. 

Set $g(t)=\text{diag}(1,t,t^2,...,t^{N-1})$, and $b_t=g(t)^{-1}bg(t)$,
$a_t=g(t)^{-1}ag(t)$. It is clear that $D(a_t,b_t)=D(a,b)$ for $t\ne 0$. 
On the other hand, since $a,b$ are triangular, $a_t,b_t$ are polynomial 
in $t$, and $a_0=0,b_0=b_{diag}$. This gives us a desired inequality, 
since the space of solutions of a system of linear equations 
cannot get bigger when the system is deformed. 
The lemma is proved. 

Now the Proposition easily follows. Indeed, conjugating 
$R(x,y)$ by the permutation $x\to x+T(x)$, we come to the situation 
$P=0,Q=0$. Thus, the proposition can be proved in the same way as 
Corollary 4.5. Proposition 4.6 is proved. 
$\square$\enddemo

\head{Appendix: 
$T$-structures and bijective 1-cocycles}\endhead

In this appendix we will introduce the notion of a $T$-structure
on an abelian group $A$, which is motivated by the definition of the map $T$
in Proposition 2.2. We will show that any group $G$ with an action 
$\rho$ on $A$ and  
a bijective 1-cocycle $\pi$ into $A$ defines such a structure, and that if $A$ 
is cyclic then the $T$-structure completely determines $G$, $\rho$, and $\pi$.  

\proclaim{Definition}
A pair $(A,T)$ of an abelian group $A$ and a bijective map $T:A\to A$
is said to be a $T$-structure if for any $x\in A$, $k\in \Z$ one has
$$
T(kx)=kT^k(x).\tag A1
$$
\endproclaim

{\bf Examples.}

1. For any abelian group, $T=id$ is a $T$-structure. 

2. Let $A=\Z$. Then there are only two $T$-structures: $T=id$ and 
$T(k)=(-1)^kk$. Indeed, let $n_{\pm}$ be such that 
$T(n_{\pm})=\pm 1$. Then $n_\pm\ne 0$ and  
$|n_\pm|T^{|n_\pm|}(\frac{n_\pm}{|n_\pm|})=\pm 1$, so $n_\pm=1$ or $n_\pm=-1$.
 If $n_+=1$, $n_-=-1$ then $T(\pm 1)=\pm 1$ and 
$T(\pm k)=|k|T^{|k|}(\pm 1)=
\pm k$ for $k>0$, so $T=id$. 
Similarly, if $n_+=-1$, $n_-=1$ then $T(k)=(-1)^kk$. 

3. Let $A$ be a ring with 1 (not necessarily commutative), 
and $c\in A$ an element such that the element 
$1+cx$ is invertible for any $x\in A$. Define $T:A\to A$ by the formula
$T(x)=x(1+cx)^{-1}$. Then $T$ is a $T$-structure on $A$ (as an additive group).
Indeed, $T$ is invertible (as $T^{-1}(y)=(1-yc)^{-1}y$), and 
it is easy to show by induction that $T^k(x)=x(1+kcx)^{-1}$, 
$T^{-k}(y)=(1-kyc)^{-1}y$, $k>0$, which proves the claim. 

4. This is a generalization of Example 2 to $A=\Z^n$. Define $A_1$ as a
subset of
$\Z^n$ consisting of integer n-tuples $(x_1,...,x_n)$ that are coprime,
i.e. with the greatest common divisor equal to  one. Denote by $A_k$ the
set of
multiples of $A_1$ by a nonnegative integer $k$.
It is easy to see that the sets $A_k$ are pairwise non-intersecting sets
whose union is the whole $\Z^n$. 
We claim that $T$-structures on $\Z^n$
are labeled by permutations of $A_1$ of order two, i.e. by the maps
$H:A_1\to A_1$, such that $H^2=id$. Indeed, starting with $H$ we first
define 
$T$ on $A_1$ simply by the formula $T(x)=H(-x),\ x\in A_1$. So defined
$T$ trivially  satisfies the relation $T(-x)=-T^{-1}(x),\ x\in A_1$. Now,
we 
extend $T$ to each $A_k$, $k$ - nonnegative, by the formula
$T(kx)=kT^k(x)$. It is easy to see that this defines a $T$-structure on
$\Z^n$. By reverting above argument and defining
$H(x)=T(-x),\ x\in A_1$ we see that all the $T$-structures on $\Z^n$ are of
this form.

The simplest properties of $T$-structures are given by the following proposition. 

\proclaim{Proposition A1}
(i) If $(A,T)$ is a $T$-structure and $A$ is finite, then $T^{|A|}=id$. 

(ii) If $(A,T)$ is a $T$-structure then $(A,T^k)$ is a $T$-structure for 
any integer $k$.

(iii) If $(A,T)$ is a $T$-structure then $(kA,\tilde T)$ is a $T$-structure
for      
any integer $k$, $\tilde T$ being the restriction of $T$ to $kA$.
\endproclaim

\demo{Proof} (i)  
Applying (A1) for $k=|A|+1$, we get 
$T((|A|+1)x)=(|A|+1)T^{|A|+1}(x)$. 
Since $|A|=0$ in $A$, the last equation yields $T(x)=T^{|A|+1}(x)$, 
thus $T^{|A|}=id$. 

(ii) Since $T^{-1}(x)=-T(-x)$, it is easy to see that $(A,T^{-1})$ is a 
$T$-structure. So it suffices to prove the statement for $k>0$. 
We do so by induction. For $k=1$ the statement known, so 
assume it is known for $k=n-1$ and let us prove it for $k=n$. 
We have $T^n(mx)=T(T^{n-1}(mx))=T(mT^{(n-1)m}(x))=
mT^m(T^{(n-1)m}(x))=mT^{mn}(x)$, as desired. 

(iii) For any $x\in A,\ k\in \Z,\ T(kx)=kT^k(x)$, so $T$ maps $kA$ to $kA$.
It follows that $(kA, \tilde T)$ is a $T$-structure.$\square$\enddemo

{\bf Remark.}  It is easy to generalize the statement (ii) of
Proposition A1. Namely, any two $T$-structures $(A,T_1),\ (A,T_2)$ on
the same group A such that $T_1T_2=T_2T_1$ give rise to a new
$T$-structure $(T, A)$ with $T=T_1T_2$.  The above observation can be
applied to Example 3. Namely, let $(A, T_1),\ (A, T_2)$ be
$T$-structures on the ring $A$, given by the formulas
$T_1(x)=x(1+c_1x)^{-1},\ T_2(x)=x(1+c_2x)^{-1}$, then
$T(x)=T_1T_2(x)=x(1+(c_1+c_2)x)^{-1}$ defines a $T$-structure.
       
Now let us explain how to construct $T$-structures by examining their
connection with bijective cocycle data.

\proclaim{Theorem A2}
Let $(G, A,\rho,\pi)$ be  a bijective cocycle
datum. Define\ \ $T~:~A~\to~A$ by the formula
$$
T(x)=\pi^{-1}(x)*x\tag A2
$$
for $x\in A$. Then, the pair $(A,T)$ is a $T$-structure.
\endproclaim

The proof of Theorem A2 follows from the Lemma stated below.

Define the product $x\circ y=\pi^{-1}(x)*y$ for $x,\ y\in A$. Then
$T(x)=x\circ x$.

\proclaim{Lemma A3} One has 
$$
(y+x)\circ z=(y\circ x)\circ (y\circ z) \tag A3
$$
for any $x,\ y\in A$.
\endproclaim

\demo{Proof of the Lemma}
One derives from the definition of a  1-cocycle that
$\pi^{-1}(y+x)=\pi^{-1}(\pi^{-1}(y)*x)\pi^{-1}(y)$.
So, $$(y+x)\circ z=\pi^{-1}(y+x)*z=
\pi^{-1}(\pi^{-1}(y)*x)*(\pi^{-1}(y)*z).$$
This implies $(y+x)\circ z=(y\circ x)\circ (y\circ z)$.$\square$\enddemo

\demo{Proof of Theorem A2} Remark 1 in Section 2.4 implies that $A$ has a 
natural structure of a nondegenerate symmetric set. The map $T$ for this 
nondegenerate symmetric set, defined in Proposition 2.2, coincides 
with the one given by (A2). Thus the 
bijectivity of $T$ follows from Proposition 2.2(b). 
 So, it remains to check that
for any integer $k \ge 0$ one has
$$T(kx)=kT^k(x),\qquad \qquad \qquad T(-kx)=-kT^{-k}(x),$$
or, in terms of "$\circ$" product, it is enough to show that
$$kx\circ x=T^k(x),  \qquad \qquad \qquad (-kx)\circ x=T^{-k}(x).\tag A4$$. 

We will prove (A4) by induction. It is clear that (A4) holds for $k=0$. 
Suppose (A4) holds for
$k=n-1$. Then, substituting $y=-nx,\ z=x$ into Lemma A2, we get
$(x-nx)\circ x=(-nT(x))\circ T(x)$, i.e. 
$T^{1-n}(x)=(-nT(x))\circ T(x)$. If we let $y=T(x)$ we obtain
$T^{-n}(y)=(-ny)\circ y$. Similarly, substituting $y=(n-1)x,\
z=x$ into Lemma A2, we get $T^n(x)=(nx)\circ x$. $\square$\enddemo

It is easy to check that examples 1-3 of $T$-structures are in fact obtained 
from bijective cocycle data by formula (A2). Indeed:

1. $T=id$ is obtained from the datum $(A,A,\text{trivial action},id)$.

2. $T(k)=(-1)^kk$ for $A=\Z$ is obtained from the datum
$(G,A,\rho,\pi)$ where $G$ is the group of affine transformations 
of $\Z$, $G=\{(a,b):x\to ax+b|a=\pm 1,b\in \Z\}$, $\rho(a,b)x=ax$, 
$\pi(a,b)=2ab+\frac{a-1}{2}$. 

3. Let $A$ be a ring with 1, and $c$ an element of $A$. 
Introduce a new operation on $A$ by $x\bullet y=
x+y+xcy$. It is easy to check that this operation is associative, and 
$0$ is the unit $e$ with respect to it. Moreover, if 
$1+cx$ is 
invertible for any $x\in A$ then 
$(A,\bullet)$ is a group: for any $x\in A$ there exists an element 
$x^{-1}=-x(1+cx)^{-1}$ such that $x\bullet x^{-1}=x^{-1}\bullet x=e$. 
Define an action $\rho: (A,\bullet)\to \End(A)$ by
$\rho(x)y=y(1+cx)^{-1}$ and 
$\pi=id:(A,\bullet)\to A$. It is easy to check 
that $((A,\bullet),A,\rho,\pi)$ is a bijective cocycle datum, 
and that the corresponding map $T$ is given by $T(x)=x(1+cx)^{-1}$.   
 
{\bf Remark.} 
Note that Remark 1 in Section 2.4 implies that in the situation of 
example 3, $A$ has a natural structure of a nondegenerate symmetric set. 
Moreover, it is easy to compute the map $S$ explicitly, and the answer is 
$$
S(x,y)=(y(1+cx+cxcy)^{-1},x(1+cy)).\tag A5
$$
Note that this formula defines the structure of a nondegenerate symmetric set 
not only on $A$ but also on any right ideal in $A$. 

In each of these examples, one notes that $T^r=id$ where
$r=|Im(\rho)|$, where $Im(\rho)$ denotes the image of $G$ in $Aut(A)$
under $\rho$.  This is a consequence of the following general
property:

\proclaim{Proposition A4}
Let $H$=$Ker(\rho)$, $K=\pi(H)$ for a bijective cocycle datum and let
$T$ be the $T$-structure obtained from the datum.  If $A/K$ is finite,
and $r$ is any integer such that $r(A/K)=0$, then $T^r=id.$
\endproclaim

\demo{Proof}
If $r$ is an integer such that $r(A/K)=0$, then $rA$ is a subset of
$K$, and hence $T^r(x)= rx \circ x = x$ since $rx \in K$, for any $x
\in A$.
$\square$\enddemo

\proclaim{Corollary A5}
In the case when $A/K$ is finite, set $r=|Im(\rho)|$, the cardinality
of the image of $G$ under $\rho$ in $Aut(A)$. Then, $T^r=id.$ In
particular, if $A$ is cyclic, then $T^{gcd(\phi(|A|),|A|)}=id$ where
$\phi$ is the Euler $\phi$-function.
\endproclaim

\demo{Proof}
Since $r=|Im(\rho)|=|A/K|$, $r(A/K)=0$ so this follows from the
previous corollary.  In the case of $A$ cyclic, we have that
$|Aut(A)|=\phi(|A|)$, $|G|=|A|$ so $|Im(\rho)|$ divides
$gcd(\phi(|A|),|A|)$.  $\square$\enddemo

Furthermore, when $T$ comes from a bijective cocycle datum as before,
then we can find a similar bijective cocycle datum that produces the
map $T^k$, which was proven to be a $T$-structure in Proposition A1,
(ii) for any integer $k$.

\proclaim{Proposition A6}
(i) Given a bijective cocycle datum $(G,A,\rho,\pi)$ and any integer
$k$, there is a unique bijective cocycle datum given by
$(G^{(k)},A,\rho^{(k)},\pi^{(k)})$ where
$\rho^{(k)}((\pi^{(k)})^{-1}(x))=\rho(\pi^{-1}(kx))$.

(ii) Let $T$ be the natural $T$-structure induced by $(G,A,\rho,\pi)$
(i.e. $T(x)=\rho(x)(x)$).  Then the natural $T$-structure induced by
$(G^{(k)},A,\rho^{(k)},\pi^{(k)})$ is given by $T^k$.
\endproclaim

\demo{Proof}
(i) Without loss of generality we may assume $G=(A,\bullet)$ and
$\pi=id$ for a binary operation $\bullet$.  We will use the notation
$x * y = \rho(x)(y)$.  Define $x \circ y = kx * y$ for $x,y \in A$ and
set $\rho^{(k)}(x)(y)=x \circ y$.  Then, define $x \ci y = y +
(-T^k(y)) \circ x$.  Now, setting $G^{(k)}=(A,\ci)$ and
$\pi^{(k)}=id$, we claim $(G^{(k)},A,\rho^{(k)},\pi^{(k)})$ is a
bijective cocycle datum.

Indeed, $(G,\ci)$ is associative:
$$(x \ci y) \ci z = x \ci (y \ci z) = z+(-T^k(z)) \circ y+((-T^k(z+(-T^k(z)))) \circ y) \circ x.$$  Furthermore, $G$ has inverses:
$$x \ci (-T^k(x))=(-T^k(x))+(-T^k(-T^k(x))) \circ x = -T^k(x) + x
\circ x = 0,$$ and $-T^k(-T^k(x))=x.$ Thus $-T^k(x)=x^{-1}$ in
$G^{(k)}$, and $G^{(k)}$ is a group.  This implies that the cocycle
condition is also satisfied: $x \ci y = y+(y^{-1}) * x$ where we take
the inverse in $G^{(k)}$. Thus, the result is a bijective cocycle
datum.

We note that this is the unique datum with the given $\rho^{(k)}$
since the operation $\ci$ is obtained directly from $\rho^{(k)}$ given
that the cocycle condition is satisfied.

(ii) It is clear that $x \circ x = T^k(x)$ since $kx * x = T^k(x)$.
$\square$\enddemo

In light of the above discussion, it makes sense to ask whether any
map $T$ comes from a bijective cocycle datum, and whether the map $T$
determines this datum. The following theorem gives a positive answer
to both questions in the case when the group $A$ is cyclic.

\proclaim{Theorem A7}
For a fixed cyclic group A, bijective cocycle data
$(G,A,\rho,\pi)$ are, up to
isomorphism, in one to one correspondence with $T$-structures on A.
\endproclaim

\demo{Proof}
Let us construct 
the group $G$, representation $\rho:G\to End(A)$ and bijective 1-cocycle
$\pi:G\to A$ from a given $T:A\to A$. 
This is done as follows. 

(1) Pick a generator $1\in A$. It defines a natural homomorphism $\Z\to A$ 
which defines a natural ring structure on $A$. Define the product
"$\circ$" on $A$ by the formula $y\circ z=zT^y(1)$, 
where we regard $A$ as a ring. Since $T^{|A|}=id$, 
$T^y$ makes sense and thus this operation is well 
defined. 

Let us check that 
$(y+z)\circ w=(y\circ z)\circ (y\circ w)$.
We have 
$$
(y+z)\circ w=wT^{y+z}(1),\ (y\circ z)\circ (y\circ w)=
wT^y(1)T^{zT^y(1)}(1).
$$
So it suffices to show that $T^z(u)=uT^{zu}(1)$ (we need the case $u=T^y(1)$), 
which was done in Proposition A1(ii). 

(2) For $y,\ z\in A$ define $y\ci z=z+(-T(z))\circ y=z+yT^{-T(z)}(1)$,
$\rho(y)z=y\circ z$. Denote by $G$ the set $A$ with the operation $\ci$. 

It is easy to see that the map $\rho:G\to End(A)$ is multiplicative. 
Indeed, we must check

$$
(y \ci z) \circ w = y \circ (z \circ w).
$$

This can be rewritten as

$$
(z \circ yT^{-T(z)}(1)) \circ (z \circ w) = y \circ (z \circ w).
$$

Now we need only to show $T^z(1) T^{-T(z)}(1) = 1$, which is true since
$T^z(1) T^{-T(z)}(1)=T^z(1) T^{-zT^z(1)}(1)=T^{z-z}(1)=1$. This also proves
that $\pi: G \rightarrow A$ satisfies the cocycle condition, since
this amounts to proving $z$ and $-T(z)$ are inverses in $G$.

This implies that $G$ is a group (the axioms of a 
group can be checked after application of $\pi$; since $\pi$ is bijective, 
the axioms must hold). Thus, we have constructed from $T$ 
a bijective cocycle datum $(G,A,\rho,\pi)$. 

Thus, we have shown that the map from bijective cocycle data to $T$-structures 
is surjective. To show that it is injective, we will show that 
if $(A,T)$ is a $T$-structure obtained from $(G,A,\rho,\pi)$ as in 
Theorem A2, and $(G',A,\rho',\pi')$ is obtained from $(A,T)$ as 
we just described, then $G'=G,\rho'=\rho,\pi'=\pi$ up to an isomorphism. 

Recall that $T$ is defined by $T(x)=\pi^{-1}(x)\circ x$. This implies that 
$T^n(x)=\pi^{-1}(nx)\circ x$. Thus $x\circ y=yT^x(1)=
y(\pi^{-1}(x)*1)=\pi^{-1}(x)*y$. Thus, 
$$
x\ci y=y+xT^{-T(y)}(1)=
y+x\pi^{-1}(-\pi^{-1}(y)*y)*1=y+\pi^{-1}(y)^{-1}*x
$$
This implies that $\pi^{-1}(x\ci y)=\pi^{-1}(x)\pi^{-1}(y)$, so 
$\pi^{-1}:G'\to G$ is a group isomorphism. This isomorphism 
clearly maps $\pi'=id$ to $\pi$. Also, 
$\rho'(x)y=\pi^{-1}(x)*y=\rho(\pi^{-1}(x))y$, so $\rho'$ corresponds to
$\rho$. The theorem is proved.  
$\square$\enddemo

\proclaim{Corollary A8} 
If $A=\Z/p\Z$, where $p$ is a prime, then any $T$-structure 
is the identity. 
\endproclaim

\demo{Proof} The group $G$ must equal $\Z/p\Z$, and its action on $A$ must 
be trivial, as there are no nontrivial actions. So $\pi$ must be a
group isomorphism. Thus, $T=id$.  Alternatively, the theorem in
conjunction with Corollary A5 immediately implies $T=id.$
$\square$\enddemo

{\bf Remark.}  This can also be proven directly as follows.  Note that $T^k(1)
\neq 0$ for all positive $k$ since $T(0)=0$.  Thus the
minimum positive $j$ such that $T^j(1)=1$ must be less than $p$ and a
factor of $p$ since $T^p(1)=1$, hence $j=1$ and $T$ is trivial since
$T(x)=xT^x(1)=x$ for all $x \in A$.

Theorem A7 raises an interesting question of classification of
$T$-structures on cyclic groups. Unfortunately, we could not find such a
classification, even when $A=\Z/n\Z$ where $n$ is a prime power.  

Note, however, by the proposition in Section 3.2, that when $A$
is cyclic $T$-structures correspond to multipermutation solutions.
In terms of the map $T$ this means that $T$ is similarly retractable
by taking $A$ modulo $Ker (\rho)$, which by the above proposition (showing $x
\circ 1 = T^x(1) = 1$ if $x \in Ker (\rho)$) is the subgroup consisting
of all elements $x \in A$ such that $T^x = id.$ This is the first step
towards such a classification.

\enddocument